\documentclass[11pt]{article}
\usepackage{a4}
\usepackage{amssymb}
\ifx\mathbb\undefined\def\mathbb{\Bbb}\fi
\makeatletter\@addtoreset{equation}{section}\makeatother
\makeatletter\@addtoreset{figure}{section}\makeatother
\makeatletter\@addtoreset{table}{section}\makeatother

\newcommand{\op}[1]{\!\!\mathop{\mbox{\rm ~{#1}}}\nolimits}
\newcommand{\fop}[1]{\!\!\mathop{\mbox{\rm\footnotesize ~{#1}}}\nolimits}
\newcommand{\scriptop}[1]{\!\!\mathop{\mbox{\rm\scriptsize ~{#1}}}\nolimits}
\newfont{\gothic}{eufm10 scaled\magstep0}
\newcommand{\R}{\mathbb R}
\newcommand{\Z}{\mathbb Z}
\newcommand{\C}{\mathbb C}
\newcommand{\Q}{\mathbb Q}
\newcommand{\proj}{\mathbb P}
\newtheorem{theorem}{Theorem}[section]
\newtheorem{proposition}[theorem]{Proposition}
\newtheorem{lemma}[theorem]{Lemma}

\newcounter{remark}[section]
\renewcommand{\theremark}{\thesection.\arabic{remark}}
\newenvironment{remark}{\refstepcounter{remark}
\par\medskip\noindent{\bf Remark~\theremark~~}}{
\unskip\nobreak\hfill\hbox{$\oslash$}\par \bigskip}

\newenvironment{proof}{\par\medskip\noindent{\bf Proof}~~}{
\unskip\nobreak\hfill\hbox{q.e.d.}\par \bigskip}    

\begin{document}
\title{Second Order Contact of Minimal Surfaces}
\author{J.J. Duistermaat}
\date{}
\maketitle

\begin{abstract}\noindent
The minimal surface equation $Q$ in the second order contact bundle 
of $\mathbb{R}^3$, modulo translations, 
is provided with a complex structure and a canonical vector\--valued holomorphic 
differential form $\Omega$ on $Q\setminus 0$. 
The minimal surfaces $M$ in $\mathbb{R}^3$ correspond 
to the complex analytic curves $C$ in $Q$, where the derivative of the 
Gauss map sends $M$ to $C$, and $M$ is equal to the real part of the 
integral of $\Omega$ over $C$. The complete minimal surfaces of finite 
topological type and with flat points at infinity correspond to the algebraic curves in $Q$.  
\end{abstract}

\noindent
{\large\bf Introduction}

\medskip\noindent
In Section \ref{n'subsec} we introduce the second order contact bundle  
modulo translations $Q$ of the minimal surface equation. 
$Q$ is a two\--dimensional 
vector bundle over the unit sphere $S$, and carries a canonical $\R ^3$\--valued 
one\--form $\omega$ which has a pole type of singularity along the zero section 
of $Q$. If $M$ is a minimal surface, then the assignment to each $x\in M$ 
of the second order contact element of $M$ at $x$ defines a mapping $n':M\to Q$, 
which can be viewed as the derivative of the Gauss map $n:M\to S$. 
The image $n'(M)$ is a two\--dimensional submanifold $C$ of $Q$ such that 
the restriction to $C$ of $\omega$ is closed. 

In Section \ref{recsubsec} we show that conversely, if $C$ is a two\--dimensional 
submanifold of $Q\setminus 0$ such that $\op{d}\!\omega |_C=0$, then the 
submanifold $M$ of $\R ^3$, which is obtained from $C$ by means of 
integration of $\omega$ over $C$, is a minimal surface in $\R ^3$ such that 
$n'(M)=C$. The integral of $\omega$ over closed curves in $C$ leads to 
a group ${\cal P}$ of periods of $M$, which is discussed in 
Section \ref{periodsubsec}. 

The next observation is that $Q$ carries a complex structure, unique up to its 
sign, such that the condition $\op{d}\!\omega |_C=0$ is equivalent to the condition 
that $C$ is a complex analytic curve in $Q$, with respect to this complex structure. 
In order to describe the complex structure, and at the same time the 
automorphism group of $Q$, we begin in Section \ref{so3subsec} by exhibiting 
$Q$ as an associated vector bundle $\op{SO}(3)\times _{\fop{SO}(2)}\R ^2$. 
In Section \ref{so3csubsec} this bundle is identified with the 
associated complex line bundle $\op{SO}(3,\,\C )\times _B\C$, 
where $B$ is a suitable Borel subgroup of $\op{SO}(3,\,\C )$. 
Because all the spaces here are complex analytic, this provides $Q$ with 
a complex structure. The subgroup $B$ is chosen in such a way that 
$\omega$ is equal to the real part of a $\C ^3$\--valued holomorphic 
$(1,\, 0)$\--form $\Omega$ 
on $Q$, equivariant for the action of $\op{SO}(3,\,\C )$ on $Q$. 

This leads to the characterization in Section \ref{Csubsec} of the $C=n'(M)$ 
as the complex analytic curves in $Q$. In Section \ref{weierstrasssubsec} we 
discuss the relation with the so\--called isotropic complex analytic curves 
in $\C ^3$, which are used in the Weierstrass type of representation formulas 
for minimal surfaces. The discussion of the structure of $Q$ is concluded 
in Section \ref{sl2subsec} with the introduction of the two\--fold covering 
$\op{SL}(2,\,\C )$ of $\op{SO}(3,\,\C )$ as a convenient computational tool. 
It leads to the identification $Q$ with the complex line bundle 
${\cal O}(4)$ over $\C\proj ^1$ and of its compactification $\overline{Q}$ 
with the fourth Hirzebruch surface $\Sigma _4$, a complex projective 
algebraic variety. It also leads to local coordinates 
in which $\Omega$ takes a relatively simple explicit form, 
cf. (\ref{uOmega}). 

In the applications, one has to pay special attention to the flat points of the 
minimal surface $M$. These correspond to certain intersection points of $C$ 
with the zero section $0_{_Q}$ of $Q$, where the one\--forms $\omega$ and 
$\Omega$ are singular. These points are discussed in Section \ref{fpsubsec}, 
whereas the other points of  $C\cap 0_{_Q}$, which correspond to 
the flat points of $M$ at infinity, are analyzed in Section \ref{inftysubsec}. 
The condition of flatness of the point at infinity is equivalent to the condition that 
the total curvature of the end of $M$ is finite. Because we allow periodicities, 
the integral of the Gaussian curvature has to be taken over the quotient 
of the end by the translation symmetry. 

In Section \ref{algsubsec} we prove that $M/{\cal P}$, the minimal surface 
modulo its periods, is of finite topological type and has only flat points 
at infinity, if and only if 
$C$ is a complex algebraic curve in the complex projective variety 
$\overline{Q}$, where $C$ has to satisfy some additional conditions 
in order to ensure that the minimal surface $M$ is smoothly immersed. 
If the Gauss map $n:M/{\cal P}\to S$ has degree $d$, then the number of 
flat points in $M/{\cal P}$ plus the number of flat points at 
infinity, each counted with multiplicity, is equal to $4d$. We obtain  
the Jorge\--Meeks formula as a consequence. In Section \ref{secsubsec} 
we study the holomorphic sections of $Q$, several of which correspond 
to known minimal surfaces. 

In Section \ref{hypsubsec}, we find that $M/{\cal P}$ 
has finite topological type, no points at infinity, and a degree two Gauss map,  
if and only if $C$ is a hyperelliptic curve of genus three, lying in the usual 
way as a twofold branched covering over $S$. The group 
${\cal P}\left(\Omega |_C\right)$ of the periods of $\Omega |_C$ is a 
lattice in $\C ^3$, and the quotient $J$, which is compact, 
is isomorphic to the Jacobian variety $\op{Jac}(C)$ of $C$.   
The integration of $\Omega$ over $C$ defines an embedding from 
$C$ onto a closed complex analytic curve in $J$.  
The isotropic complex curve $\Gamma$ in $\C ^3$ such that $M=\op{Re}\Gamma$ 
as in Section \ref{weierstrasssubsec}, is equal to the pull\--back under 
the projection $\C ^3\to J$ of the image of $C$ in $J$. Therefore $\Gamma$ 
is a closed and smooth one\--dimensional complex  analytic submanifold of 
$\C ^3$. A result of Pirola \cite{pirola} implies that the 
hyperelliptic curves $C$, for which the period group 
${\cal P}$ of the minimal surface $M$ is a discrete subgroup of $\R ^3$, 
are dense. In an open dense subset they form a countable union of 
families of hyperelliptic curves which depend on essentially 
five real parameters. 

There is an enormous literature on minimal surfaces, among  which 
the beautiful surveys of Nitsche and those of 
Dierkes, Hildebrandt, K\"uster and Wohlrab. 
We will use DHKW \cite{dhkw} as our main reference. 

\section{Second order Contact}
\label{order2sec}
In this section we introduce the real two\--dimensional vector 
bundle $Q$ over the unit sphere $S$ such that the combination of the 
Gauss map and its derivative sends each minimal surface $M$ in 
$\R ^3$ to a surface $C$ in $Q$. Moreover, there is a canonically 
defined $\R ^3$\--valued one\--form $\omega$ on $Q$ such that 
the minimal surface $M$ is reconstructed from $C$ by means of integration 
of $\omega$ over curves in $C$.  
\subsection{The Derivative of the Gauss Map}
\label{n'subsec}
Let $M$ be an oriented two\--dimensional smooth submanifold of $\R ^3$ 
and $S=\{ s\in\R ^3\mid\| s\| =1\}$ the sphere in $\R ^3$ with 
radius equal to one and center at the origin. The mapping 
$n:M\to S$ which assigns to $x\in M$ the oriented normal $n(x)$ to 
the tangent space $\op{T}_xM$, 
is called the {\em Gauss map} of $M$. The tangent map 
$\op{T}_xn$ of $n$ at $x$ is a linear mapping from $\op{T}_xM$ 
to $\op{T}_{n(x)}S$. The translation $\tau _{x-n(x)}$ from $n(x)$ to 
$x$ leads to an identification of $\op{T}_{n(x)}S$ with 
$\op{T}_xM$. Let $n'(x)=\op{T}_{x}n\circ\tau _{x-n(x)}$ 
denote the linear mapping from $\op{T}_{n(x)}S$ to $\op{T}_{n(x)}S$ 
which is induced by $\op{T}_xn$. 

The linear endomorphism 
$n'(x)$ of $\op{T}_{n(x)}S$ is symmetric with respect to 
the restriction to $\op{T}_{n(x)}S$ of the standard inner 
product of $\R ^3$. In fact, by means of a translation 
we can arrange that a given special point of $M$ is at the 
origin, and by rotation in $\R ^3$ we can subsequently arrange that 
$n(0)=e_3$, the vertical standard basis vector in $\R ^3$. 
Then, near $0$, $M$ can be written as the graph 
of a smooth function $f$ of two variables, the third 
coordinate as a function of the first two ones. If we identify 
$\op{T}_{e_3}S$ with $\R ^2\simeq\R ^2\times\{ 0\}$, then 
a short calculation shows that the matrix of $n'(0)$ is equal to 
minus the Hessian matrix $f''(0)$ of $f$ at $0$. Because $f''(0)$ 
describes the second order contact of $M$ at $0$ with its tangent 
plane, we have in general that $n'(x)$ represents the second 
order contact of $M$ at $x$ with $\op{T}_xM$. 

The surface $M$ is a {\em minimal surface} if and only if, 
for every $x\in M$, the trace of $n'(x)$ 
is equal to zero, cf. DHKW \cite[(24) on p. 16 and p. 53]{dhkw}. 
For every $s\in S$, let 
\begin{equation}
Q_s:=\{ q\in\op{Lin}\left(\op{T}_sS,\,\op{T}_sS\right)\mid 
q^*=q,\;\op{trace}q=0\} 
\label{Ws}
\end{equation}
be the space of all traceless symmetric linear mappings from 
$\op{T}_sS$ to itself. $Q_s$ is a two\--dimensional real vector 
space. The $Q_s$, $s\in S$, form a smooth rank two real vector 
bundle $Q$ over $S$. If $M$ is a minimal surface then, for 
every $x\in M$, the element $n'(x)\in Q_{n(x)}$ represents 
the second order contact of $M$ at $x$. This defines a smooth 
mapping $n':M\to Q$ such that the Gauss mapping $n$ is equal to 
$\pi\circ n'$, if $\pi :Q\to S$ denotes the projection which 
assigns to $q\in Q_s$ the base point $s$. 

Let $M$ be a minimal surface in $\R ^3$. 
A traceless symmetric two by two matrix is of the form 
$q=\left(\begin{array}{cc}a&b\\b&-a\end{array}\right)$ in which 
$a,\, b\in\R$. Its determinant is equal to $-a^2-b^2$ and it 
follows that if $q\neq 0$, then $q$ is invertible. 
This implies that, for every $x\in M$ 
the {\em Gauss curvature} 
\begin{equation}
K(x)=\op{det}n'(x)\quad \left(\mbox{\rm equal to}\; -a^2-b^2
\;\mbox{\rm when}\; n(x)=e_3\right) 
\label{Kab}
\end{equation}of $M$ at $x$ 
is nonpositive. Moreover, if $x$ is an {\em umbilic point of $M$}, 
i.e. $K(x)=0$, then $x$ is a {\em flat point of $M$}, 
i.e. $n'(x)=0$, which is equivalent to the condition that $M$ osculates 
its tangent plane at the point $x$. 
In other words, the following conditions i)---v) are equivalent. 
\begin{itemize}
\item[i)] $x$ is not a flat point of $M$. 
\item[ii)] $n'(x)\neq 0$. 
\item[iii)] $n'(x)$ is invertible. 
\item[iv)] The Gauss map 
$n$ is a diffeomorphism from some open neighborhood $M_0$ of $x$ in $M$ 
onto an open neighborhood $S_0$ of $n(x)$ in $S$. 
\item[v)] There is an open neighborhood $M_0$ of $x$ in $M$ 
such that the restriction to $M_0$ of $n'$ is a smooth 
embedding from $M_0$ to a smooth local section $C_0=n'\left( M_0\right)$ 
of the vector bundle $Q$ over $S$.
\end{itemize}
Here the equivalence between iv) and v) follows from the fact 
that $n=\pi\circ n'$. Indeed, this implies that if $n$ is a diffeomorphism, then 
$n'$ is an embedding. Conversely if $n'$ is a diffeomorphism from 
$M_0$ onto a smooth submanifold $C_0$ of $W$, then $C_0$ is 
a smooth local section of $\pi :Q\to S$ $\Longleftrightarrow$ the 
restriction to $C_0$ of $\pi$ is a diffeomorphism from $C_0$ onto 
an open subset $S_0$ of $S$ $\Longleftrightarrow$ $\pi\circ n'$ 
is a diffeomorphism from $M_0$ onto $S_0$. 

\subsection{Reconstruction of the Minimal Surface by Integration}
\label{recsubsec}
If $q\in Q_s\setminus\{ 0\}$, then $q^{-1}:\op{T}_sS\to\op{T}_sS$ 
can be viewed as a linear mapping from $\op{T}_sS$ to $\R ^3$, 
if we identify the target tangent plane as a linear subspace of $\R ^3$. 
In this way we obtain an $\R ^3$\--valued one\--form 
(= differential form of degree one) on $Q\setminus 0$, which 
we denote by $\omega =q^{-1}\,\op{d}\! s$. 
Here $0=0_{_{Q}}$ denotes the zero section of $Q$. More precisely, 
if $q\in Q_s$, then $\omega$ is the linear mapping from 
$\op{T}_qQ$ to $\R ^3$ which is equal to the linear 
mapping $\op{T}_q\pi$ from $\op{T}_qQ$ to $\op{T}_sS$, 
followed by the linear mapping $q^{-1}$ from $\op{T}_sS$ to 
$\R ^3$. Note that $\omega$ 
has a pole\--type of singularity at $0_{_Q}$. 

Suppose that $M$ is a connected immersed minimal surface and 
$n'$ is a diffeomorphism 
from $M$ onto a smooth submanifold $C$ of $Q$. 
Let $\iota :M\to\R ^3$ denote the identity on $M$, viewed 
as a mapping from $M$ to $\R ^3$. Then, for every $x\in M$ and 
$v\in\op{T}_xM$, we have 
\[
\left((n')^*\omega\right) _x(v)=\omega _{n'(x)}\left(\op{T}_xn'(v)\right) 
=q^{-1}\circ\op{T}_{n'(x)}\pi\circ\op{T}_xn'(v)
=q^{-1}\circ\op{T}_xn(v)=v,
\]
where $q=n'(x)=\op{T}_xn$ is viewed as a linear mapping from 
$\op{T}_xM\simeq\op{T}_{n(x)}S$ to $\R ^3$. 
Here we have applied the chain rule to $\pi\circ n'=n$ 
in the third identity. This proves that 
\begin{equation}
(n')^*\omega =\op{d}\!\iota\quad\mbox{\rm on}\quad M.
\label{omegaiota}
\end{equation} 
In turn this means that $\omega |_C=\op{d}\! f$ if 
$f$ is equal to the $\R ^3$\--valued function 
$\left((n')^{-1}\right) ^*\iota =(n')^{-1}$ on $C$. 
It follows that $M$ can be reconstructed from $C$ 
in the following way. Choose, for any base point $q_0\in C$, a corresponding 
base point $x_0\in\R ^3$ as $(n')^{-1}\left( q_0\right)$. 
Then, for every $q_1\in V$, the point of $M$ corresponding to $q_1$ 
is given by  
\begin{equation}
x_1=\phi\left( q_1\right) 
:=(n')^{-1}\left( q_1\right) =x_0+\int_{q_0}^{q_1}
\,  q^{-1}\,\op{d}\! s,
\label{q-1}
\end{equation}
in which the integral of the $\R ^3$\--valued one\--form 
$\omega =\underline{q}^{-1}\,\op{d}\! s$ 
is taken along any smooth curve $\gamma$ in $C$ which runs from 
$q_0$ to $q_1$. Note that the integral does not depend on the 
choice of $\gamma$. 

The fact that $\omega |_C$ is exact, is locally equivalent to the 
condition that $(\op{d}\!\omega )|_C=
\op{d}\! \left(\omega |_C\right) =0$, which means that 
the restriction to $C$ of $\op{d}\!\omega$ is equal to zero. 
This condition means that, for every $q\in C$ and every 
pair $v,\, w\in\op{T}_qC$ of tangent vectors to $C$, 
we have that $(\op{d}\!\omega )_q(v,\, w)=0$, or that 
$\op{T}_qC$ is an isotropic linear subspace of $\op{T}_qQ$ 
with respect to the $\R ^3$\--valued antisymmetric 
bilinear form $(\op{d}\!\omega )_q$. 
Because it is clear that $\op{d}\left( q^{-1}\,\op{d}\! s\right) 
=\op{d}\left( q^{-1}\right)\wedge\op{d}\! s$ is far from zero 
as a two\--form on the four\--dimensional manifold $Q\setminus 0$, 
the equation $(\op{d}\!\omega )|_C=0$ is strong restriction 
on the two\--dimensional real linear subspace $\op{T}_qC$ 
of $\op{T}_qQ$.  

\subsection{Surfaces in the Bundle and their Periods}
\label{periodsubsec}
Suppose that $C$ is any smooth two\--dimensional immersed 
submanifold of $Q\setminus 0$ such that $(\op{d}\!\omega )|_C=0$, 
which means that $\omega |_C$ is closed. 
Also assume that $C$ is transversal to the fibers of the 
projection $\pi :Q\to S$, which means that 
the $\R ^3$\--valued linear form $q^{-1}\,\op{d}\! s$ 
is injective on $\op{T}_qC$. Then, 
after a choice of base points $q_0\in C$, $x_0\in\R ^3$, 
(\ref{q-1}) defines an embedding $\phi$ 
from each simply connected open subset $C_0$ of $C$ 
onto a smooth two\--dimensional submanifold $M_0$ of $\R ^3$. 
For each $q\in C$ we have, with the notation $s=\pi (q)$, 
$x=\phi (q)$, that 
\[
\op{T}_xM_0=\op{image}\op{T}_q\phi =q^{-1}\left(\op{T}_sS\right) 
=\op{T}_sS,
\]
where $\op{T}_{\phi (q)}M_0$ and $\op{T}_sS$ both are regarded 
as two\--dimensional linear subspaces of $\R ^3$. It follows 
that $s$ is orthogonal to $\op{T}_xM_0$, and we have a unique orientation of 
$M_0$ such that $s=n(x)$, the normal of $M_0$ at the point $x$. 
Furthermore the equation $\op{d}\! x=q^{-1}\,\op{d}\! s$ 
implies that $\op{d}\! s=q\,\op{d}\! x$. Therefore 
the inverse of $\phi$ is equal to $n'$. The fact that $n'$ maps 
$M_0$ into $Q$, where the fiber $Q_s$ is equal to 
the space of traceless symmetric linear mappings from 
$\op{T}_sS$ to $\op{T}_sS$, finally implies that $M_0$ is 
a minimal surface. In this way {\em the mapping $n'$ establishes 
a local equivalence between the minimal surfaces in $\R ^3$ 
with nonvanishing curvature, 
modulo translations, and the two\--dimensional smooth submanifolds 
$C$ of $Q\setminus 0$ such that $\omega |_C$ is closed and 
$C$ is transversal to the fibers of the projection $\pi :Q\to S$}. 

Let $[\omega]$ denote the de Rham cohomology class of $\omega$, 
which is an element of $\op{H}_{\scriptop{de Rham}}^1(C)\otimes\R ^3$. 
Globally, the equation (\ref{q-1}) in general defines a 
multi\--valued immersion 
$\phi$ from $C$ to $\R ^3$, where the multi\--valuedness 
is caused by the fact that for every closed loop $\gamma$ in $C$, 
meaning that $q_1=q_0$, 
we obtain that $x_1-x_0=\langle [\gamma ],\, [\omega ]\rangle$ 
need not be equal to zero. The $\langle [\gamma ],\, [\omega]\rangle$ , 
where $[\gamma ]$ denotes the homology class of the closed loop 
$\gamma$ in the image of $\op{H}_1(C,\,\Z )$ in $\op{H}_1(C,\,\R )$, 
are called the {\em periods} 
of $\omega$. They form an additive subgroup ${\cal P}$ of $\R ^3$, 
generated by the vectors $\langle\left[\gamma _j\right] ,\, [\omega]\rangle$
if the $[\gamma _j]$ generate the image of $\op{H}_1(C,\, \Z )$ 
in $\op{H}_1(C,\,\R )$. 
The multi\--valuedness of $\phi$ means that $\phi$ assigns 
to every $q$ a coset of the form $x+{\cal P}$ in $\R ^3$. 
This implies that the image $M$ of $C$ is {\em ${\cal P}$\--periodic} 
in the sense that $M+p=M$ for every $p\in {\cal P}$. 
Conversely, if $M+p=M$, then $n'(x+p)=n'(x)$ for every 
$x\in M$, and it follows that the set of periods of $M$ is equal to 
${\cal P}$. 

In Section \ref{complexsec} we will show that $Q$ has a complex 
structure, unique up to sign, 
such that if $q\in Q\setminus 0$ then  
a two\--dimensional real linear subspace $P$  
of $\op{T}_qQ$ is isotropic with respect to 
$(\op{d}\!\omega )_q$ if and only if $P$ is equal to a 
one\--dimensional complex linear subspace of the two\--dimensional 
complex vector space $\op{T}_qQ$. This implies that a surface $C$ 
is of the form $n'(M)$ for a minimal surface $M$ with nonvanishing 
curvature if and only if 
$C$ is a complex analytic curve in the two\--dimensional 
complex analytic manifold $Q\setminus 0$ which is transversal to the fibers. 

\begin{remark}
Any second order partial differential equation for surfaces in $\R ^3$ 
can be identified with a 7\--dimensional hypersurface in the 8\--dimensional 
second order contact bundle of $\R ^3$. If the equation 
is translation invariant, then we can pass to its quotient 
by the translation group $\R ^3$, which is a 4\--dimensional manifold $Q$. 
If the equation is elliptic, then there is an almost complex structure 
$J$ on $Q$ such that the solution surfaces, modulo translations, 
are locally in a bijective correspondence with the complex analytic 
curves in $(Q,\, J)$. For the minimal surface equation the 
almost complex structure $J$ is integrable. I actually arrived at the 
description of the complex structure in $Q$ from the application 
of the theory of second order contact structure to the minimal 
surface equation. However, in order to emphasize the especially 
nice features of the minimal surface equation and not to 
burden the presentation with the generalities about second
order contact structures, I have chosen to present the direct 
description of the complex structure of Section \ref{complexsec}.  
\end{remark}

\begin{remark}
The considerations in this section have been predominantly of a local 
nature. Avoiding the flat points, the mapping $n'$ is an immersion from 
$M$ to $Q\setminus 0$ which need not be injective, let alone that it is an 
embedding. One cause for non\--injectivity could be the occurrence 
of nonzero periods, which can be remedied by passing to the quotient 
$M/{\cal P}$. However, there can also be self\--intersections of the 
image $C=n'(M)$ which are not caused by periodicities; at such  
multiple points $C$ has a singularity of multiple point type. 
Finally the 
mapping $n':M/{\cal P}\to Q\setminus 0$ need not be proper, 
meaning that one can have a non\--converging sequence of points in 
$M$ such that the image points converge to a point in $Q\setminus 0$. 

In the other direction, the immersion $(n')^{-1}:C\to\R ^3$ can have 
self\--intersections which generically occur along curves. In this 
paper we will not discuss the problem of avoiding such 
self\--intersections of $M$ or $M/{\cal P}$. 
\end{remark}

\section{The Complex Structure}
\label{complexsec}
In this section, a complex structure is introduced on $Q$ by 
identifying $Q$ with an associated complex line bundle 
over $\op{SO}(3,\,\C)\times _B\C$, where $B$ is the 
group of all $b\in\op{SO}(3,\,\C )$ such that 
$b\left(\C\,\left( e_1+\op{i}e_2\right)\right) =
\C\,\left( e_1+\op{i}e_2\right)$. Here $b\in B$ 
acts on $\C$ by means of multiplication by $\chi (b)$, 
in which $\chi :B\to\C\setminus\{ 0\}$ is the unique 
character (= Lie group homomorphism) such that 
$\chi (b)=\op{e}^{\fop{i}2\phi}$ if $b$ is equal to the 
rotation abouth the vertical axis through the angle $\phi$. 
In this description the base space $S$ of $Q$, the unit sphere 
in $\R ^3$, is identified with $\op{SO}(3,\,\C )/B$, 
which in turn is identified with the quadric $N$ in 
the complex projective plane, defined by the 
homogeneous quadratic equation $\langle z,\, z\rangle =0$, $z\in\C ^3$. 

The $\R ^3$\--valued one\--form $\omega$ on $Q\setminus 0$ 
is equal to the real part of an $\op{SO}(3,\,\C )$\--equivariant 
$\C ^3$\--valued holomorphic $(1,\, 0)$\--form 
$\Omega$. For a real two\--dimensional submanifold $C$ of 
$Q\setminus 0$, we have that $\omega |_C$ is closed if 
and only if $C$ is a complex analytic curve in the 
complex two\--dimensional manifold $Q\setminus 0$. 
If $M$ is the minimal surface in $\R ^3$ which is obtained by 
means of integration of $\omega |_C$, and $\Gamma$ 
is the complex analytic curve in $\C ^3$ which is obtained 
by means of integration of $\Omega |_C$, then 
$M=\op{Re}\Gamma$, and $\Gamma$ is the isotropic 
complex anaytic curve corresponding to $M$ in the 
Weiestrass type of representation formulas. 

The computations simplify considerably if one replaces the 
group $\op{SO}(3,\,\C )$ by the group $\op{SL}(2,\,\C )$. 
The adjoint representation (= action on the Lie algebra 
$\mbox{\gothic sl}(2,\,\C )$ of $\op{SL}(2,\,\C )$ 
by means of conjugation) leads to a two\--fold covering 
$\op{Ad}:\op{SL}(2,\,\C )\to 
\op{SO}(3,\,\C )$, in which we use a suitable identification 
of $\mbox{\gothic sl}(2,\,\C )$ with $\C ^3$. 
Under this identification, $B$ corresponds to 
the group $L$ of all lower triangular 
matrices in $\op{SL}(2,\,\C )$, and 
$Q\simeq\op{S0}(3,\,\C )\times _B\C$ to 
the associated line bundle $\op{SL}(2,\,\C )\times _L\C$, 
where $\left(\begin{array}{cc}z&0\\v&1/z\end{array}\right)\in L$ 
acts on $\C$ by means of multiplication by $z^4$. 
In this description, the sphere $S\simeq N$ is identified with 
the complex projective line ${\C}{\proj}^1$, of which 
$\op{SL}(2,\,\C )/\{\pm 1\}$ is the automorphism group. 
It follows that the automorphism group of $Q$ is equal to 
$\left(\op{SL}(2,\,\C )/\{\pm 1\}\right)\times\left(\C\setminus\{ 0\}\right)  
\simeq\op{SO}(3,\,\C )\times\left(\C\setminus\{ 0\}\right)$.  
Here $c\in\C\setminus\{ 0\}$ acts on $Q$ by multiplication with 
$c$ in each fiber. 

\subsection{The Rotation Group}
\label{so3subsec}
The vector bundle $Q$ is homogeneous with respect to the group $\op{SO}(3)$ 
of rotations in $\R ^3$, in the sense that the 
mapping $(g,\, q)\mapsto \left( g\left( e_3\right) 
,\, g\, q\, g^{-1}\right)$ 
is surjective from $\op{SO}(3)\times Q_{e_3}$ to $Q$. 
We have that $(g,\, q)$ and $(g', q')$ have the same image if 
and only if there exists a rotation $r$ about the vertical 
axis, such that $g'=g\, r^{-1}$ and $q'=r\, q\, r^{-1}$. 
Therefore, if we denote by $\op{SO}(2)$ the group of all 
rotations about the vertical axis, then $Q$ is isomorphic 
to the space $\op{SO}(3)\times _{\fop{SO}(2)}Q_{e_3}$ 
of $\op{SO}(2)$\--orbits in $\op{SO}(3)\times Q_{e_3}$, 
where $r\in\op{SO}(2)$ acts on $\op{SO}(3)\times Q_{e_3}$ 
by sending $(g,\, q)$ to $\left( g\, r^{-1},\, r\, q\, r^{-1}\right)$. 
The projection $(g,\, q)\mapsto g$ factorizes to a 
projection from $\op{SO}(3)\times _{\fop{SO}(2)}Q_{e_3}\simeq Q$ 
onto $\op{SO}(3)/\op{SO}(2)\simeq S$, which corresponds to the 
projection $\pi :Q\to S$. Here $\op{SO}(3)/\op{SO}(2)$ is 
identified with $S$ by means of the mapping 
$g\mapsto g\left( e_3\right)$, which is surjective 
from $\op{SO}(3)$ onto $S$ and has the $\op{SO}(2)$\--orbits  
as its fibers. The element $h\in\op{SO}(3)$ acts on $Q$ 
by sending $(s,\, q)$ to $\left( h\, s,\, h\, q\, h^{-1}\right)$, 
which corresponds to sending $(g,\, q)$ to $(h\, g,\, q)$. 
For this reason this action is called the {\em left action} of $h$ 
on $Q$ and denoted by $\op{L}_h$. 
The vector bundle $Q\simeq\op{SO}(3)\times _{\fop{SO}(2)}Q_{e_3}$ is 
also called the {\em associated $\op{SO}(3)$\--bundle 
for the representation $r\mapsto \left( q\mapsto r\, q\, r^{-1}\right)$ 
of $\op{SO}(2)$ on the vector space $Q_{e_3}$}, cf. 
\cite[Sec. 2.4]{dk}. 

Note that if 
$r=\left(\begin{array}{cc}\cos\phi&-\sin\phi
\\ \sin\phi&\cos\phi\end{array}\right)$
is the rotation in the horizontal plane through the angle 
$\phi$, and 
$q=\left(\begin{array}{cc}a&b\\b&-a\end{array}\right)$, 
then $g\, q\, g^{-1}=\left(\begin{array}{cc}a'&b'
\\b'&-a'\end{array}\right)$, in which 
$a'=(\cos 2\phi )\, a-(\sin 2\phi )\, b$ and 
$b'=(\sin 2\phi )\, a+(\cos 2\phi )\, b$.
In other words, the vector $(a',\, b')$ 
is obtained from the vector $(a,\, b)$  
by applying the rotation in $\R ^2$ through the angle $2\phi$. 

The $\R ^3$\--valued differential form $\omega =q^{-1}\,\op{d}\! s$ is 
equivariant for the left action of $\op{SO}(3)$ on $Q$, in the sense that 
if $h\in\op{SO}(3)$, then 
\begin{eqnarray*}
\left(\left(\op{L}_h\right)^*\omega\right) _{(s,\, q)}(\delta s,\,\delta q)
&=&\omega _{\left( h\, s,\, h\, q\, h^{-1}\right)} 
\left( h\,\delta s,\,\dots\right) \\
&=&\left( h\, q\, h^{-1}\right) ^{-1}
\, h\,\delta s=h\, q^{-1}\,\delta s
=h\,\omega _{(s,\, q)}(\delta s,\,\delta q).
\end{eqnarray*}

\subsection{The Complex Rotation Group}
\label{so3csubsec}
At $s=e_3$ and  
a nonzero element $q$ of $Q_{e_3}$, $\omega _{\left( e_3,\, q\right)}$ sends 
$(\delta s,\,\delta q)$ to 
the vector $y=q^{-1}\,\delta s\in\R ^3$, such that $y_3=0$. 
If $q=\left(\begin{array}{cc}a&b\\b&-a\end{array}\right)$ then 
$q^{-1}=\frac{1}{a^2+b^2}\,\left(\begin{array}{cc}a&b\\b&-a\end{array}\right)$ 
and therefore 
\begin{eqnarray*}
y_1&=&\frac{1}{a^2+b^2}\left( a\,\delta s_1+b\,\delta s_2\right) 
=\frac{1}{a^2+b^2}\,\op{Re}\left[ 
(a-\op{i}b)\,\left(\delta s_1+\op{i}\delta s_2\right)\right] 
=\op{Re}\,\frac{\delta s_1+\op{i}\delta s_2}{a+\op{i}b},\\
y_2&=&\frac{1}{a^2+b^2}\left( b\,\delta s_1-a\,\delta s_2\right) 
=\frac{1}{a^2+b^2}\,\op{Re}\left[\op{i}\,  
(a-\op{i}b)\,\left(\delta s_1+\op{i}\delta s_2\right)\right] 
=\op{Re}\,\op{i}\,\frac{\delta s_1+\op{i}\delta s_2}{a+\op{i}b}.
\end{eqnarray*}

If we identify $Q_{e_3}$ with $\C$ by identifying $q=\left(\begin{array}{cc}a&b\\b&-a\end{array}\right)$ with 
$a+\op{i}b$, and if we write $\delta s=\delta g\left( e_3\right)$ with 
$\delta g\in\mbox{\gothic so}(3)$, then we obtain that  
\begin{equation}
\widetilde{\omega}_{(g,\, q)}(\delta g,\,\delta q)
=\op{Re}\left[\widetilde{\Omega} _{(g,\, q)}(\delta g,\,\delta q)\right] , 
\label{omegaOmega}
\end{equation}
in which 
\begin{equation}
\widetilde{\Omega} _{(g,\, q)}(\delta g,\,\delta q)
=\frac{1}{q}\,
\langle g^{-1}\circ\delta g\left( e_3\right),\, e_1+\op{i}e_2\rangle\cdot 
g\left( e_1+\op{i}e_2\right) .   
\label{Omega}
\end{equation}
 
Here $\widetilde{\Omega}$ is an equivariant $\C ^3$\--valued 
one\--form on $\op{SO}(3)\times\C\setminus\{ 0\}$, 
which moreover is equal to the pull\--back of a one\--form 
$\Omega$ on 
$Q\setminus 0
\simeq\op{SO}(3)\times _{\fop{SO}(2)}\C\setminus\{ 0\}$ 
by means of the projection from $\op{SO}(3)\times\C\setminus\{ 0\}$ 
onto $Q\setminus 0$. The rotation $r$ about the vertical axis 
through the angle $\phi$ acts on $\C$ by means of multiplication 
by $\op{e}^{\fop{i}2\phi}$, because $q\mapsto r\, q\, r^{-1}$ 
corresponds to applying $r^2$ to the first column of $q\in Q_{e_3}$. 

If $V$ is a complex vector space, then a $V$\--valued one\--form 
$\Theta$ on a complex manifold $P$ is called 
{\em a $(1,\, 0)$\--form} if, for every $p\in P$, $\Theta _p$ is a 
complex linear mapping from $\op{T}_pP$ to $V$. If moreover 
in local holomorphic coordinates the coefficients of $\Theta _p$ 
depend holomorphically on $p$, then $\Theta$ is called a 
$V$\--valued holomorphic $(1,\, 0)$\--form on $P$. 
If we extend the standard inner product of 
$\R ^3$ to the corresponding complex bilinear form on $\C ^3\times\C ^3$, 
then (\ref{Omega}) extends to a $\C ^3$\--valued 
holomorphic $(1,\, 0)$\--form 
on the complex analytic manifold $\op{SO}(3,\,\C )\times\C\setminus\{ 0\}$, 
which we also denote by $\widetilde{\Omega}$. 
The idea is to introduce a closed Lie subgroup $B$ of 
$\op{SO}(3,\,\C )$ with the following properties: 
\begin{itemize}
\item[i)] $\op{SO}(2)\subset B$. 
\item[ii)] The 
injections $\op{SO}(3)\to\op{SO}(3,\,\C )$ and $\op{SO}(2)\to B$ 
induce an isomorphism from $Q\simeq\op{SO}(3)\times _{\fop{SO}(2)}\C$ 
onto $\op{SO}(3,\,\C )\times _B\C$. 
\item[iii)] $\widetilde{\Omega}$ 
is equal to the pull\--back of a differential form $\Omega$ 
on $\op{SO}(3,\,\C )\times _B\C$ by means of the projection 
$\psi :\op{SO}(3,\,\C )\times\C\to\op{SO}(3,\,\C )\times _B\C$. 
\end{itemize}
As a consequence, $\op{SO}(3,\,\C )\times _B\C$ will be a holomorphic 
complex line bundle over the complex one\--dimensional complex analytic 
manifold $\op{SO}(3,\,\C )/B$, which must be isomorphic to the 
complex projective line ${\C}{\proj}^1$ because it is diffeomorphic 
to the sphere $\op{SO}(3)/\op{SO}(2)\simeq S$. Furthermore, 
$\Omega$ then is an $\op{SO}(3,\,\C )$\--equivariant 
$\C ^3$\--valued holomorphic $(1,\, 0)$\--form 
on $\op{SO}(3,\,\C )\times _B\C$, such that $\op{Re}\Omega =\omega$  
in the identification $Q\stackrel{\sim}{\to}\op{SO}(3,\,\C )\times _B\C$. 

The condition that $\widetilde{\Omega}=\psi ^*\Omega$ for some 
$\Omega$ implies that $\widetilde{\Omega}$ is equal to zero in the 
direction of the $B$\--orbit, which means in view of (\ref{Omega}) 
that 
\[
\langle X\left( e_3\right) ,\, e_1+\op{i}e_2\rangle =0
\]
for every element $X$ of the Lie algebra $\mbox{\gothic b}$ of $B$. 
This means that $X$ must be of the form 
\begin{equation}
X=\left(\begin{array}{ccc}0&-a&-b\\
a&0&-\fop{i}b\\
b&\fop{i}b&0\end{array}\right) ,\quad a,\, b\in\C .
\label{Xab}
\end{equation}
which are precisely the $X\in\mbox{\gothic so}(3,\,\C )$ 
such that $X\left( e_1+\op{i}e_2\right)$ is equal to a 
complex multiple of $e_1+\op{i}e_2$, where for  
(\ref{Xab}) the factor is equal to $-\op{i}a$. 
For this reason, we define $B$ as the set of all 
$b\in\op{SO}(3,\,\C )$ such that 
\begin{equation}
b\left(\C\,\left( e_1+\op{i}e_2\right)\right) =\C\,\left( e_1+\op{i}e_2\right) .
\label{Bdef}
\end{equation}
$B$ is a closed complex Lie subgroup of $\op{SO}(3,\,\C )$, 
with Lie algebra $\mbox{\gothic b}$ 
equal to the set of $X$ as in (\ref{Xab}). 

Let $\widetilde{N}$ denote the {\em isotropic cone in $\C ^3$}, 
which consists of the $z\in\C ^3$, such that 
$\langle z,\, z\rangle =0$. Let $N$ denote the corresponding 
quadric in the complex projective plane, which consists of 
the $\C\, z$ such that $z\in\widetilde{N}\setminus\{ 0\}$. 
Because $e_1+\op{i}e_2\in\widetilde{N}\setminus\{ 0\}$ and $\op{SO}(3,\,\C )$ 
leaves $\widetilde{N}\setminus\{ 0\}$ invariant, we have a mapping 
$g\mapsto g\left(\C\,\left( e_1+\op{i}e_2\right)\right)$ 
from $\op{SO}(3,\,\C )$ to $N$, which in view of 
the definition (\ref{Bdef}) of $B$ induces an injective 
mapping from $\op{SO}(3,\,\C )/B$ to $N$. 

In order to find a natural diffeomorphism $\sigma$ from $N$ onto $S$, 
we write $z\in\C ^3$, $z\neq 0$, as $z=x+\op{i}y$ with 
$x,\, y\in\R ^3$. Then $z\in\widetilde{N}$ if and only if 
$\langle x,\, x\rangle =\langle y,\, y\rangle$ and 
$\langle x,\, y\rangle =0$, which implies that 
$\| x\wedge y\| =\| x\| ^2$, where $x\wedge y\in\R ^3$ denotes the 
exterior product of $x$ and $y$. If $a,\, b\in\R$, then 
$(a+\op{i}b)\, (x+\op{i}y)=a\, x -b\, y+\op{i}\, (b\,x -a\, y)$ 
and $(a\, x-b\, y)\wedge (b\, x+a\, y)=\left( a^2+b^2\right)\, x\wedge y$. 
For every $z\in\widetilde{N}\setminus\{ 0\}$, we write  
\begin{equation}
\sigma (z)=\langle x,\, x\rangle ^{-1}\, x\wedge y. 
\label{sigma}
\end{equation}
Then $\sigma (z)\in S$ and $\sigma (c\, z)=\sigma (z)$ for 
every $c\in\C\setminus\{ 0\}$, and therefore $\sigma :\widetilde{N}\to S$ 
induces a mapping from $N$ to $S$, which we also denote by 
$\sigma$. This mapping is 
$\op{SO}(3)$\--equivariant, and therefore surjective because 
$\op{SO}(3)$ acts transitively on $S$. The mapping is injective.  
Proof: if $\sigma (x+\op{i}y)=\sigma (x'+\op{i}y')$ and
$x+\op{i}\, y,\, x'+\op{i}y'\in\widetilde{N}$, then 
both $\| x\| ^{-1}\, x,\, \| y\| ^{-1}\, y$ and 
$\|x'\| ^{-1}\, x',\,\| y'\| ^{-1}\, y'$ form an orthonormal basis  
of the same two\-dimensional linear subspace of $\R ^3$ with the 
same orientation. Therefore these are obtained from each other by means 
of a rotation in this plane, which in turn implies that 
$x'+\op{i}y'\in\C (x+\op{i}y)$. The conclusion is that 
$\sigma :N\to S$ is a diffeomorphism. 

As a consequence, $\op{SO}(3)$   
acts transitively on $N$, which implies that $\op{SO}(3,\,\C )$ 
acts transitively on $N$. We conclude that the mapping 
$g\mapsto g\left(\C\,\left( e_1+\op{i}e_2\right)\right)$ induces 
an isomorphism from $\op{SO}(3,\,\C )/B$ onto $N\simeq S$. 
Because $\sigma\left(\C\,\left( e_1+\op{i}e_2\right)\right) =e_3$ 
and the action of $g\in\op{SO}(3)$ commutes with $\sigma$, 
it follows also that $B\cap\op{SO}(3)=\op{SO}(2)$. 

Like $\op{O}(3)$, the group $\op{O}(3,\,\C )$ has two 
connected components, corresponding to the sign of the 
determinant, cf. Chevalley \cite[p. 16]{c}. Therefore $\op{SO}(3,\,\C )$ 
is connected. Because $\op{SO}(3,\,\C )/B\simeq N\simeq S$ is 
simply connected, it follows that $B$ is connected as well. 
The definition of $B$ in (\ref{Bdef}) implies that, for 
every $b\in B$, $b\left( e_1+\op{i}e_2\right) =
\rho (b)\,\left( e_1+\op{i}e_2\right)$, in which 
$\rho$ is a character of $B$, a Lie group homomorphism 
from $B$ to the muliplicative group of the nonzero 
complex numbers. Because $B$ is connected, the character 
$\rho$ is determined by the infinitesimal character 
$\rho ':\mbox{\gothic b}\to\C$, where $\rho '(X)=\, -\op{i}a$ 
if $X$ is as in (\ref{Xab}). In the definition of 
$Q\simeq\op{SO}(3)\times _{\fop{SO}(2)}\C$ the action of 
the rotation $r$ about the vertical axis 
through the angle $\phi$ on $\C$ is by means of multiplication 
by $\op{e}^{\fop{i}2\phi}$. Therefore the inclusions 
$\op{SO}(3)\to\op{SO}(3,\,\C )$ and $\op{SO}(2)\to B$ 
lead to a well\--defined mapping from 
$\op{SO}(3)\times _{\fop{SO}(2)}\C$ to 
$\op{SO}(3,\,\C )\times _B\C$, if and only if we let 
act $b\in B$ act on $\C$ by means of multiplication by 
$\rho (b)^{-2}$. The isomorphism $\op{SO}(3)/\op{SO}(2)\to 
\op{SO}(3,\,\C )/B$ then implies that the mapping 
from $\op{SO}(3)\times _{\fop{SO}(2)}\C$ to 
$\op{SO}(3,\,\C )\times _B\C$ is an isomorphism of 
complex line bundles. 

The last property of $\widetilde{\Omega}$ which 
has to verified in order that $\widetilde{\Omega}=\psi ^*\Omega$ for some 
$\Omega$, is that $\widetilde{\Omega}$ is invariant under the 
right action $\op{R}_b$ of $b\in B$. We have in view of (\ref{Omega}) that 
\begin{eqnarray*}
{\op{R}_b}^*\widetilde{\Omega}_{(g,\, q)}(\delta g,\,\delta q)
&=&\left(\rho (b)^{-2}\, q\right) ^{-1}\,
\langle b\circ g^{-1}\circ\delta g\circ b^{-1}\left( e_3\right) 
,\, e_1+\op{i}e_2\rangle
\cdot g\circ b^{-1}\left( e_1+\op{i}e_2\right)\\
&=&\frac{1}{q}\,\rho (b)^2\,\langle g^{-1}\circ\delta g\circ b^{-1}
\left( e_3\right) , 
\, b^{-1}\left( e_1+\op{i}e_2\right)\rangle 
\cdot g\circ b^{-1}\left( e_1+\op{i}e_2\right) \\
&=&\frac{1}{q}\,\langle g^{-1}\circ\delta g\circ b^{-1}\left( e_3\right) ,
\, e_1+\op{i}e_2\rangle 
\cdot g\left( e_1+\op{i}e_2\right) \\
&=&\frac{1}{q}\,\langle g^{-1}\circ\delta g\circ b^{-1}\left( e_3\right) , 
\, e_1+\op{i}e_2\rangle 
\cdot g\left( e_1+\op{i}e_2\right)
= \widetilde{\Omega}_{\left( g,\, q\right)}
\left(\delta g,\delta q\right) ,
\end{eqnarray*}
which proves the desired invariance. Here we have used 
in the first equality that $\left( g\circ b^{-1}\right) ^{-1}= 
b\circ g^{-1}$ and that 
$\delta\left( g\circ b^{-1}\right) =\delta g\circ b^{-1}$. 
In the third equality we have used that 
$b^{-1}\,\left( e_1+\op{i}e_2\right) =
\rho (b)^{-1}\,\left( e_1+\op{i}e_2\right)$.  
In the fourth 
equality we have used that $b\left( e_3\right)$ is equal to 
$e_3$ plus a multiple of $e_1+\op{i}e_2$, because 
the elements of $\mbox{\gothic b}$ map $e_3$ into 
$\C\,\left( e_1+\op{i}e_2\right)$, cf. (\ref{Xab}). 
Note also that $\langle g^{-1}\circ\delta g\left( e_1+\op{i}e_2\right) , 
\, e_1+\op{i}e_2\rangle =0$ because $g^{-1}\circ\delta g$ is 
antisymmetric. 
\begin{remark}
We have, for every $g\in\op{SO}(3,\,\C )$, that 
$g\,\left( e_1+\op{i}e_2\right)\in\widetilde{N}$, which implies in view of 
(\ref{Omega}) that $\widetilde{\Omega}$, and therefore $\Omega$ as well, 
takes its values in the isotropic cone $\widetilde{N}$ rather 
than in $\C ^3$. Let $C$ be a complex analytic 
curve in $Q\setminus 0$. Let $\Gamma$ be the complex analytic 
curve in $\C ^3$ which is obtained by integration of 
$\Omega$ over curves in $C$, as in (\ref{q-1}) with $\omega =q^{-1}\,\op{d}\! s$ 
replaced by $\Omega$. Then $\Gamma$ is an {\em isotropic curve} 
in the sense that, for every $z\in\Gamma$, we have that $\op{T}_z\Gamma\in N$. 
(Note however that this does not imply that  
$\widetilde{z}-z\in\widetilde{N}$ for every $z,\,\widetilde{z}\in\Gamma$.)
\label{conerem}
\end{remark}
\begin{remark}
The group $\op{SO}(3,\,\C )$ acts transitively on $Q\setminus 0\simeq 
\op{SO}(3,\,\C )\times _B\C\setminus\{ 0\}$ and 
$\C\left( e_1+\op{i}\, e_2\right)$ is the unique one\--dimensional complex linear 
subspace of $\C ^3$ which is fixed by $B$. It follows that every 
$\op{SO}(3,\,\C )$\--equivariant $\C ^3$\--valued $(1,\, 0)$\--form on 
$Q\setminus 0$ is equal to a constant multiple of $\Omega$. I owe this 
observation to Erik van den Ban. 
\label{erikrem}
\end{remark}

\subsection{Complex Analytic Curves in Q}
\label{Csubsec}
We now determine the real two\--dimensional (immersed) submanifolds 
$C$ of  $Q$ such that $\op{d}\!\omega |_C=0$, where $\omega =\op{Re}\Omega$. 

Because of the $\op{SO}(3)$\--equivariance of $\omega$, 
we only need to investigate the problem which real two\--dimensional 
linear subspaces $P$ of $\op{T}_qQ$ are isotropic with respect to 
$(\op{d}\!\omega )_q$ at points $q\in Q_{e_3}$, which corresponds to 
$g=I$ in the identification of  $q$ with the $B$\--orbit of $(g,\, q)$ in 
$\op{SO}(3,\,\C )\times\C$. It follows from (\ref{Omega}) that 
$\op{d}\!\widetilde{\Omega}$ is nowhere zero on $\op{SO}(3,\,\C )\times\C$, 
which implies that $\op{d}\!\Omega$ is nowhere zero on $Q$, because 
$\op{d}\!\widetilde{\Omega}=\op{d}\!\left(\psi ^*\Omega\right)  
=\psi ^*(\op{d}\!\Omega ) $, if $\psi$ denotes the projection from 
$\op{SO}(3,\,\C )\times\C$ onto $Q\simeq\op{SO}(3,\,\C )\times _B\C$. 
Therefore, $(\op{d}\!\Omega )_q=\Theta\, (1,\,\op{i},\, 0)$, in which 
$\Theta$ is a nonzero antisymmetric complex bilinear form on $\op{T}_qQ$. 

If $a,\, b\in\op{T}_qQ$, then 
\[
0=\omega _q(a,\, b)=\op{Re}\Omega _q(a,\, b)
=\op{Re}\left(\Theta (a,\, b)(1,\,\op{i},\, 0)\right) 
\]
$\Longleftrightarrow\left(\op{Re}\Theta (a,\, b)=0\& 
\op{Re}\left(\Theta (a,\, b)\,\op{i}\right) =0\right) 
\Longleftrightarrow\Theta (a,\, b)=0$. \\
Because $\op{T}_qQ$ is a two\--dimensional 
complex vector space and $\Theta$ is a nonzero antisymmetric complex bilinear 
$\C$\--valued form on $\op{T}_qQ$, we have that $\Theta (a,\, b)=0$ 
if and only if 
$a$ and $b$ are linearly dependent over $\C$. 

\begin{proposition}
If $C$ is a real two\--dimensional submanifold of $Q$, 
then $\omega |_C$ is closed if and only if $C$ is a complex analytic curve 
in $Q$. 

If $M$ is 
a smoothly immersed minimal surface in $\R ^3$, provided with 
the complex structure which makes the Gauss map $n:M\to S$ 
complex analytic, then the mapping $n':M\to Q$ is also complex analytic.
\label{ccprop}
\end{proposition} 

\begin{proof}
Let $P$ be a real two\--dimensional linear subspace of $\op{T}_qQ$. 
Then $\op{d}\!\omega _p$ is equal to zero on $P\times P$ if and only 
if $\Theta |_{P\times P}=0$, which is certainly the case if $P$ is a 
one\--dimensional complex\--linear subspace of  $\op{T}_qQ$. 
Conversely, if  $\Theta _{P\times P}=0$ and $a\in P$, $a\neq 0$, then 
every $b\in P$ is a complex multiple of $b$, which implies that 
$P$ is a complex\--linear subspace of $\op{T}_qQ$. Therefore, 
$\op{d}\!\omega _p$ is equal to zero on $P\times P$ if and only 
if $P$ is a one\--dimensional complex\--linear subspace of $\op{T}_qQ$. 

The second statement follows from the ``only if'' part of the first statement. 
\end{proof}

\begin{remark}
The (almost) complex structure on $Q$ with the above property 
is unique up to its sign. Indeed, 
suppose that $J$ is a complex structure in $\op{T}_qQ$, which 
means that $J$ is a real linear mapping from $\op{T}_qQ$ to itself, 
such that $J^2=\, -I$. Suppose moreover that every 
one\--dimensional complex linear subspace $P$ of $\op{T}_qQ$ 
is also a complex\--linear subspace with respect to $J$, 
which means that $J(P)=P$. This implies that for every 
nonzero $a\in\op{T}_qQ$ there exist unique real numbers 
$\alpha =\alpha (a)$ and $\beta =\beta (a)$, depending 
smoothly on $a$, such that $J(a)=\alpha (a)\, a+\beta (a)\,\op{i}a 
=\gamma (a)\, a$, where we have written 
$\gamma (a)=\alpha (a)+\beta (a)\op{i}$. 
Replacing $a$ by $a+b$ such that $a$ and $b$ are linearly 
independent over $\C$,  
we obtain from $J(a+b)=J(a)+J(b)$ that 
$\gamma (a+b)\, (a+b)=\gamma (a)\, a+\gamma (b)\, b$, 
hence $\gamma (a)=\gamma (a+b)=\gamma (b)$. Because these 
$(a,\, b)$ are dense in $\op{T}_qQ\times\op{T}_qQ$, 
it follows that $\gamma$ is a constant. But then $J^2=-I$ 
is equal to multiplication by means of the constant $\gamma ^2$, 
which implies that $\gamma =\pm\op{i}$. This proves that  
the complex structure $J$ on $\op{T}_qQ$, such that 
the set of the $(\op{d}\!\omega )_q$\--isotropic two\--dimensional 
real linear subspaces of $\op{T}_qQ$ is equal to the set 
of the complex one\--dimensional linear subspaces with 
respect to $J$, is uniquely determined up to its sign. 
\end{remark}

\subsection{The Weierstrass Representation}
\label{weierstrasssubsec}
Let $M$ be an oriented two\--dimensional real submanifold of 
$\R ^3$ and  $x\in M$. Then $\op{T}_xM$  has a unique complex structure 
$J_x$ such that $J_x$ is antisymmetric with respect 
to the Euclidean inner product in $\op{T}_xM$ and such that, 
for every nonzero $v\in\op{T}_xM$, the pair $v,\, J_x(v)$ is 
positively oriented. 
(In view of 
${J_x}^2=-I$, the antisymmetry of $J_x$ is equivalent to the condition that $J_x$ is 
an orthogonal linear transformation.) 

Recall the identity mapping $\iota :M\to\R ^3$, viewed as an $\R ^3$\--valued 
smooth function on $M$, which had been used 
in (\ref{omegaiota}). For each $x\in M$, we view the $\R ^3$\--valued one\--form 
$a:=\op{d}\!\iota _x$ on $\op{T}_xM$ as a linear mapping from $\op{T}_xM$ to $\R ^3$. 
There is a unique complex\--linear mapping 
$\Phi _x:\op{T}_xM\to\C ^3$ such that, for every 
$v\in\op{T}_xM$, $\op{d}\!\iota _x(v)=\op{Re}\Phi _x(v)$. 
Indeed, if we write $b(v)=\op{Im}\Phi _x(v)$, then 
$v\mapsto a(v)+\op{i}b(v)$ is complex\--linear 
if and only if 
$\op{i}\left( a(v)+\op{i}b(v)\right)$ is equal to  
$a\left( J_x(v)\right) +\op{i}b\left( J_x(v)\right)$, 
which in turn is equivalent to $a(v)=b\left( J_x(v)\right)$ and 
$-b(v)=a\left( J_x(v\right)$. The second equation 
implies the first in view of ${J_x}^2=\, -I$, and we arrive at the 
conclusion that 
\begin{equation}
\Phi _x:=\op{d}\!\iota _x-\op{i}\,\op{d}\!\iota _x\circ J_x
\label{Phi}
\end{equation}
is the unique complex\--linear mapping 
$\Phi _x:\op{T}_xM\to\C ^3$ such that 
$\op{d}\!\iota _x=\op{Re}\Phi _x$. In other words, there is a unique 
$\C ^3$\--valued $(1,\, 0)$\--form $\Phi$ on $M$, such that 
\begin{equation}
\op{d}\!\iota =\op{Re}\Phi .
\label{diotaPhi}
\end{equation}

\begin{remark}
It follows from (\ref{Phi}) that, for each $v\in\op{T}_xM$, 
the real and imaginary part of $\Phi _x(v)$ have the same length 
and are orthogonal to each other, which is equivalent 
to the statement that $\langle\Phi _x(v),\,\Phi _x(v)\rangle =0$.  
In other words, {\em $\Phi$ takes its values in the isotropic cone 
$\widetilde{N}$} rather than in $\C ^3$. In turn this implies that, 
for each $x\in M$, $\Phi _x\left(\op{T}_xM\right)\in N$, where $N$ is  
the quadric in the complex projective plane defined by the equation 
$\langle z,\, z\rangle =0$. 
\label{PhiNrem}
\end{remark}

If $M$ is a minimal surface with nonzero curvature, then 
$M$ has a unique complex structure $J$ such that the Gauss map 
$n:M\to S$ is complex analytic, where we have provided $S$ 
with the complex structure via the diffeomorphism $\sigma :N\to S$  
defined by (\ref{sigma}). Note that $q=n'(x)$ reverses the 
orientation if we identify $\op{T}_xM$ and $\op{T}_{n(x)}S$ 
by means of a translation, which means that 
$n(x)=\, -\langle v,\, v\rangle ^{-1}\, v\wedge J_x(v)$, 
if $v\in\op{T}_xM$ and $v\neq 0$. The definitions 
(\ref{Phi}) and (\ref{sigma}) therefore imply that the 
Gauss map can be expressed in terms of $\Phi$ by means of 
\begin{equation}
n(x)=\sigma\left(\Phi _x\left(\op{T}_xM\right)\right) ,\quad x\in M,
\label{nsigmaPhi}
\end{equation}
cf. DHKW \cite[(22) on p. 94]{dhkw}. 

It follows from (\ref{diotaPhi}) and (\ref{omegaiota}) that 
$
\op{Re}\Phi =\op{d}\!\iota =(n')^*\omega 
=(n')^*\op{Re}\Omega =\op{Re}\left( (n')^*\Omega\right) ,
$
which in view of the uniqueness of the complexification of 
$\op{d}\!\iota$ implies that 
\begin{equation}
\Phi =(n')^*\Omega .
\label{PhiOmega}
\end{equation}
Note that both $\Phi$ and  $\Omega$ 
are $\widetilde{N}$\--valued, cf. Remark \ref{conerem}.  
Because $n'$ is a holomorphic mapping from $M$ to $Q$ and 
$\Phi$ is holomorphic on $Q\setminus 0$, 
we recover in this way the classical fact that {\em $M$ is a minimal surface if 
and only if the $\C ^3$\--valued $(1,\, 0)$\--form $\Phi$ on $M$ is 
holomorphic}, cf. DHKW \cite[Sec. 2.6 and 3.1]{dhkw}. The equation 
(\ref{PhiOmega}) implies that {\em all the one\--forms 
$\Phi$ on the different minimal manifolds $M$ are pull\--backs 
of one and the same canonical $\C ^3$\--valued holomorphic 
$(1,\, 0)$\--form $\Omega$ on the fixed space $Q$}. Only the 
mapping $n'$, with which $\Omega$ is pulled back, varies with 
the minimal surface. 

Any $(1,\, 0)$\--form on a one\--dimensional complex analytic 
manifold is holomorphic if and only if it is closed. 
Suppose that, for any base point $x_0\in M$, we have chosen 
a point $z_0\in\C ^3$ 
such that $\op{Re}z_0=x_0$. 
Then we define, for every $x\in M$, the point $z=z(x)\in\C ^3$ by 
\begin{equation}
z=z_0+\int_{x_0}^{x}\,\Phi ,
\label{zx}
\end{equation}
in which the integral of the $\C ^3$\--valued one\--form 
$\Phi$ is taken along any smooth curve $\gamma$ in $M$ which runs from 
$x_0$ to $x$. Note that the point $z$ is only unique modulo 
periods $\langle [\gamma],\, [\Phi ]\rangle$, where 
$[\gamma ]$ denotes the homology class, in the image of 
$\op{H}_1(M,\,\Z )$ in $\op{H}_1(M,\,\R )$, 
of a closed loop $\gamma$ in $M$ and $[\Phi ]$ denotes the 
de Rham cohomolgy class of $\Phi$ in 
$\op{H}^1_{\scriptop{de Rham}}(M)\otimes\C ^3$. 
It follows from (\ref{diotaPhi}) and $\op{Re}z_0=x_0$ that 
\begin{equation}
\op{Re}z(x)=x,\quad x\in M,
\label{Rez=x}
\end{equation}
which in turn implies that the periods of (\ref{zx}) are purely 
imaginary. The periods of (\ref{zx}) form an additive subgroup of $\op{i}\R ^3$, 
generated by the $\langle [\gamma _j],\, [\Phi ]\rangle$, 
if the $\left[\gamma _j\right]$ generate the 
image of $\op{H}_1(M,\,\Z )$ in $\op{H}_1(M,\,\R )$. 

The $z(x)$, $x\in M$, form a complex analytic curve 
$\Gamma$ in $\C ^3$, such that $M=\op{Re}\Gamma$. 
Here $\op{Re}:\C ^3\to\R ^3$ is the projection 
which sends $z=x+\op{i}y\in\C ^3$, with $x,\, y\in\R ^3$, 
to its real part $x=\op{Re}z$. We have 
$\op{T}_{z(x)}\Gamma =\Phi _x\left(\op{T}_xM\right)$.  
It follows from Remark \ref{PhiNrem} that 
the curve $\Gamma$ is {\em isotropic} in the sense 
that, for every $z\in\Gamma$, we have that $\op{T}_z\Gamma\in N$. 

Conversely, for 
every isotropic complex analytic curve $\Gamma$ in $\C ^3$, the 
real two\--dimensional manifold $M=\op{Re}\Gamma$ is a minimal 
surface. This leads to an identification between minimal surfaces 
in $\R ^3$ and isotropic complex analytic curves in $\C ^3$. 
The various {\em representation formulas of Weierstrass type}, as 
discussed in DHKW \cite[Sec. 3.1 and 3.3]{dhkw}, consist of 
constructions of isotropic complex analytic curves in $\C ^3$. 
These are obtained by the integration of a suitable 
complex analytic function 
$\phi :D\to\C ^3$, where $D$ is an open subset of $\C$, 
such that, for every $\zeta\in D$, $\phi (\zeta )\neq 0$ and 
$\langle\phi (\zeta ),\,\phi (\zeta )\rangle =0$. 

Let $C$ be a complex analytic curve in $Q\setminus 0$. Let $\Gamma$ be the 
complex analytic curve in $\C ^3$ which is obtained by 
integration of $\Omega |_C$, and let $M=\op{Re}\Gamma$ be the minimal surface in 
$\R ^3$, which is obtained by integration of  $\omega |_C=\op{Re}\Omega|_C$. 
That is, $\Gamma$ is 
defined by (\ref{q-1}), with 
$\omega =q^{-1}\,\op{d}\! s$ replaced by $\Omega$. 
It then follows from $C=n'(M)$ and (\ref{PhiOmega}) that 
$\Gamma$ is equal to the complex curve in $\C ^3$ which is obtained from 
the minimal surface $M$ in $\R ^3$ by integration of  the $\C ^3$\--valued 
holomorphic $(1,\, 0)$\--form $\Phi$ on $M$.

\subsection{SL(2,\, C)}
\label{sl2subsec}
For computations, the group $\op{SL}(2,\,\C )$, of all complex $2\times 2$\--matrices 
with determinant equal to one, is easier to work with than the complex rotation group 
$\op{SO}(3,\,\C )$. 
The element $g\in\op{SL}(2,\,\C )$ acts on 
its Lie algebra $\mbox{\gothic sl}(2,\,\C )$ of 
traceless complex $2\times 2$\--matrices by means of the 
conjugation $\op{Ad}g :X\mapsto g\, X\, g^{-1}$. 
 
If we use the identification 
\begin{equation}
\Xi (x)=\left(\begin{array}{cc}\fop{i}x_3&\fop{i}x_1-x_2
\\\fop{i}x_1+x_2&-\op{i}x_3\end{array}\right) ,\quad x\in\C ^3. 
\label{Xx}
\end{equation}
of $\C ^3$ with $\mbox{\gothic sl}(2,\,\C )$, then 
\begin{equation}
\op{trace}(\Xi (x)\, \Xi (y))=\, -2\langle x,\, y\rangle ,\quad x,\, y\in\C ^3.
\label{Xxy}
\end{equation}
Because the trace is invariant under conjugation, it follows that 
the adjoint representation (= the action on the Lie algebra by means 
of conjugation) defines a Lie group homomorphism 
$\op{Ad}:\op{SL}(2,\,\C )\to \op{SO}(3,\,\C )$.

Because $\op{ad}=\op{T}_I\op{Ad}$ is injective, the adjoint 
representation is a covering map from $\op{SL}(2,\,\C )$ 
onto a subgroup of $\op{SO}(3,\,\C )$ of the same dimension 
as $\op{SL}(2,\,\C )$. Because both $\op{SO}(3,\,\C )$ and 
$\op{SL}(2,\,\C )$ are complex 3\--dimensional and 
$\op{SO}(3,\,\C )$ is connected, it follows that 
$\op{Ad}\left(\op{SL}(2,\,\C )\right) =\op{SO}(3,\,\C )$. 
Because the kernel of the adjoint representation of 
$\op{SL}(2,\,\C )$ is equal to $\{\pm I\}$, the adjoint 
representation defines a two\--fold covering from $\op{SL}(2,\,\C )$ 
onto $\op{SO}(3,\,\C )$, or equivalently an isomorphism 
\begin{equation}
\op{Ad}:\op{SL}(2,\,\C )/\{\pm I\}\stackrel{\sim}{\to}\op{SO}(3,\,\C ).
\label{Adiso}
\end{equation}

\begin{lemma}
$\op{Ad}^{-1}(B)$ is equal to the group $L$ of all 
lower triangular matrices in $\op{SL}(2,\,\C )$.
The homomorphisms $\op{Ad}:\op{SL}(2,\,\C )\to\op{SO}(3,\,\C )$ 
and $\op{Ad}:L\to B$ induce an isomorphism 
\begin{equation}
\op{SL}(2,\,\C )\times _L\C\stackrel{\sim}{\to}
\op{SO}(3,\,\C )\times _B\C\simeq Q,
\label{LB}
\end{equation}
if we let act the element 
$g=\left(\begin{array}{cc}a&0\\c&1/a\end{array}\right)$ 
on $\C$ by means of multiplication by $a^4$. 
The pull\--back of $\Omega$ under (\ref{LB}) is an equivariant 
$\C ^3\simeq\mbox{\gothic sl}(2,\,\C )$\--valued $(1,\, 0)$\--form 
on $\op{SL}(2,\,\C )\times _L\C\setminus\{ 0\}$, which we also denote by $\Omega$. Its pull\--back to $\op{SL}(2,\,\C )\times \C\setminus\{ 0\}$ is given by 
\begin{equation}
\widetilde{\Omega}_{(g,\, q)}(\delta g,\,\delta q)=\frac{\fop{i}}{q}\,\op{trace}
\left( g^{-1}\,\delta g\, Y\right)\cdot g\, Y\, g^{-1}. 
\label{Omegasl2}
\end{equation}
\label{sl2so3lem}
\end{lemma}

\begin{proof}
We have $Y:=\Xi \left( e_1+\op{i}e_2\right) =\left(\begin{array}{cc}0&0\\2\fop{i}&0\end{array}\right)$. 
Furthermore, if $g\in\op{SL}(2,\,\C )$ then $g\, Y\, g^{-1}=\rho\, Y$ for 
some $\rho\in\C$ if and only if 
$g=\left(\begin{array}{cc}a&0\\c&1/a\end{array}\right)$ for some 
$a\in\C\setminus\{ 0\}$ and $b\in\C$, in which case $\rho =a^{-2}$. 
This proves the first statement. 
 
In the definition of  $\op{SO}(3,\,\C )\times _B\C$, the element 
$b\in B$ acts on $\C$ by means of multiplication by $\rho ^{-2}$, 
if $b\,\left( e_1+\op{i}e_2\right) =\rho\,\left( e_1+\op{i}e_2\right)$. 
This implies the second statement. 

For the proof of (\ref{Omegasl2}), we 
recall the definition (\ref{Omega}) of 
$\Omega$ on $\op{SO}(3,\,\C )\times _B\C\setminus\{ 0\}$. 
The infinitesimal action of 
$X=g^{-1}\,\delta g=\left(\begin{array}{cc}a&b\\c&-a\end{array}\right)$ on 
$\Xi\left( e_3\right) =\left(\begin{array}{cc}\fop{i}&0\\0&-\fop{i}\end{array}\right)$ 
is equal to the commutator 
$
X\circ\Xi\left( e_3\right)-\Xi\left( e_3\right)\circ X =\left(\begin{array}{cc}0&-2\fop{i}b\\2\fop{i}c&0\end{array}\right) 
$
of  $X$ and $\Xi\left( e_3\right)$.
Its product with $Y$ is equal to 
$\left(\begin{array}{cc}4b&0\\0&0\end{array}\right)$, of which the trace 
divided by $-2$, corresponding to the inner product in (\ref{Omega}), is equal to $-2b$. 
On the other hand the trace of 
$X\, Y$ is equal to $2\op{i}b$, which multiplied by $\op{i}$ yields $-2b$, 
which completes the proof of (\ref{Omegasl2}). 
\end{proof}

\begin{remark}
In view of (\ref{Xxy}), the isotropic cone $\widetilde{N}$ of the $x\in\C ^3$ 
such that $\langle x,\, x\rangle =0$ corresponds to the set of 
$X\in\mbox{\gothic sl}(2,\,\C )$ such that $\op{trace}\left( X^2\right) =0$, 
which are precisely the nilpotent elements in $\mbox{\gothic sl}(2,\,\C )$. 
Note that the element $Y$ in (\ref{Omegasl2}) is nilpotent. In other words, 
the $\mbox{\gothic sl}(2,\,\C )$\--valued one\--form $\Omega$ on 
$\op{SL}(2,\,\C )\times _L\C$ actually takes its values in 
the {\em nilpotent cone} in $\mbox{\gothic sl}(2,\,\C )$. 
\end{remark}

The element $g\in\op{SL}(2,\,\C )$ acts on the complex projective 
line ${\C}{\proj}^1$ by sending the one\--dimensional 
complex linear subspace $l$ of $\C ^2$ to $g(l)$. 
If $e_2$ denotes the second standard basis vector in 
$\C ^2$ and $l=\C\, e_2$, then $g(l)=l$ if and only if $g\in L$. 
Because the action of $\op{SL}(2,\,\C )$ 
on ${\C}{\proj}^1$ is transitive, 
the mapping $g\mapsto g(l)$ induces a 
isomorphism from $\op{SL}(2,\,\C)/L$ 
onto ${\C}{\proj}^1$. This leads to an identification of  
\begin{itemize}
\item[i)] the unit sphere 
$S\simeq\op{SO}(3)/\op{SO}(2)$ in $\R ^3$, 
\item[ii)] the quadric 
$N\simeq\op{SO}(3,\,\C )/B$ in the complex projective plane, and  
\item[iii)]
the complex projective line ${\C}{\proj}^1\simeq\op{SL}(2,\,\C )/L$ 
\end{itemize}
with each other. 

The standard projective 
coordinate in the neighborhood ${\C}{\proj}^1\setminus\{\C\, e_1\}$ 
of $\C\, e_2$ in ${\C}{\proj}^1$, is obtained by sending $u\in\C$ to 
the coset $g\, B\in\op{SL}(2,\,\C )/L$ of $g$, where 
$g=\left(\begin{array}{cc}1&u\\0&1\end{array}\right)$. 
The mapping which sends $(u,\, q)\in\C ^2$ to the $L$\--orbit 
of  $(g,\, q)\in\op{SL}(2,\,\C )\times\C$ is the corresponding 
trivialization of the complex line bundle $Q$. 

We obtain that 
$g^{-1}\,\delta g=\left(\begin{array}{cc}0&\delta u\\0&0\end{array}\right)$, 
the product of which with $Y$ has trace equal to $2\op{i}\delta u$.  Furthermore 
\begin{equation}
g\, Y\, g^{-1}=\left(\begin{array}{cc}2\fop{i}u&-2\fop{i}u^2\\
2\fop{i}&-2\fop{i}u\end{array}\right) 
=\Xi (x), \quad x=\left( 1-u^2,\,\op{i}\,\left( 1+u^2\right) ,\, 2u\right) . 
\label{AdgY}
\end{equation}
It follows that, in the $(u,\, q)$\--coordinates on $Q$, 
\begin{equation}
\Omega =\frac{-2}{q}\,\op{d}\! u\,\left( 1-u^2,\, 
\op{i}\,\left( 1+u^2\right) ,\, 2u\right) . 
\label{uOmega}
\end{equation}
Recall that $\omega$ is equal to the real part of (\ref{uOmega}). 

\begin{remark}
Suppose that the complex curve $C$ is locally equal to the graph 
of a complex analytic function $q=q(u)$, in which $u$ is the 
coordinate on ${\C}{\proj}^1$ which has been used in (\ref{uOmega}). 
If we write ${\cal F}(u)=\, -2/q(u)$, then (\ref{uOmega})  
coincides with the {\em representation formula of Weierstrass},  
which is given in DHKW \cite[p.113]{dhkw}.  

However, the formula (\ref{uOmega}) is more explicit, because of the 
interpretation of the variable $u$ as an analytic coordinate of the unit sphere, 
a coordinate for the image of the Gauss map, whereas $q$ represents 
the second order contact of  the minimal surface $M$, if   
$C=n'(M)$ and $M$ is obtained from $C$ by means of integration 
of $\omega |_C=\Omega |_C$. 
\label{FWrem}
\end{remark}

\begin{remark}
It follows from (\ref{AdgY}) that 
\[
\sigma\left(\C\, g\, Y\, g^{-1}\right)  =\left( 1+{u_1}^2+{u_2}^2\right) ^{-1}\,
\left( -2u_1,\, -2u_2,\, 1-{u_1}^2-{u_2}^2\right) ,\quad u=u_1+\op{i}u_2, 
\]
in which $\sigma$ denotes the isomorphism from the quadric $N$ in the 
complex projective plane to the unit sphere $S$ as defined in (\ref{sigma}).  
The right hand side is the formula for the stereographic projection from the 
plane to the unit sphere minus the south pole $-e_3$, but with a reversed orientation. 
The reversal of the 
orientation corresponds to the reversal of orientation of 
the derivative of the Gauss map. 
\label{stereographrem}
\end{remark}

The standard projective coordinate in the coordinate neighborhood of 
the missing point  is obtained by sending $v\in\C$ to 
the coset $h\, B\in\op{SL}(2,\,\C )/L$ of $g$, where 
$h=\left(\begin{array}{cc}0&1\\-1&v\end{array}\right)$. 
If  $g=\left(\begin{array}{cc}1&u\\0&1\end{array}\right)$ 
and $k=\left(\begin{array}{cc}a&0\\b&1/a\end{array}\right)$, 
then  $h\, k=g$ if and only if  $a=1/u$ and $b=u$.  Because 
$(g,\, q)$ and $\left( h,\, a^4\, q\right)$ are in the same $L$\--orbit in 
$\op{SL}(2,\,\C )\times \C$, they define the same point in 
$\op{SL}(2,\,\C )\times _L\C\simeq Q$, which means that 
\begin{equation}
(u,\, q)\mapsto \left( u^{-1},\, u^{-4}\, q\right) 
\label{uv}
\end{equation} 
is the corresponding coordinate transformation (= retrivialization) of $Q$. 
The complex line bundle $Q$ can be constructed by glueing two copies 
of $\C ^2$ together on $\left(\C\setminus\{ 0\}\right)\times\C$ by means 
of the mapping (\ref{uv}).   

As a consequence, the global holomorphic sections of $Q$ are of 
the form $q=f_0(u)$ where $f_0$ is a polynomial of degree fours, 
cf. Section \ref{secsubsec}. In other words, the Chern number 
or degree of the holomorphic line bundle $Q$ over $\C\proj ^1$ 
is equal to four, or $Q\simeq {\cal O}(4)$, 
cf. Griffiths and Harris \cite[pp. 144, 145]{gh}. 

\begin{remark}
The description with $\op{SL}(2,\,\C )$ leads to an easy determination 
of the automorphism group of the complex line bundle $B$. 
Let  $\Phi$ be an automorphism of the complex line bundle 
$Q$, in the sense that $\Phi$ is a complex analytic diffeomorphism 
on $Q$ which maps each fiber to another fiber by means of a complex linear mapping. 
Viewing the fibers as the points of the base space ${\C}{\proj}^1$, this 
leads to a complex analytic diffeomorphism $\Psi$ on ${\C}{\proj}^1$. 
It is known that each such $\Psi$ is equal to the action 
of an element $g$ of $\op{SL}(2,\,\C)$, where $g$ is uniquely determined 
up to it sign. See Griffiths and Harris \cite[p. 64]{gh}. 
Let $\Phi _0$ be equal to the composition of  
$\Phi$ and the left action of $g^{-1}$ on $Q\simeq\op{SL}(2,\,\C )\times _L\C$. 
Then $\Phi _0$ maps each fiber of $Q$ to itself by means of multiplication 
by a nonzero complex number, which depends in a complex analytic fashion 
on the base point in ${\C}{\proj}^1$. Because ${\C}{\proj}^1$ is compact, it 
follows from the maximum principle that this function is a constant. 

We conclude that the automorphism group of the complex line bundle $Q$ 
is equal to 
$\left(\op{SL}(2,\,\C )/{\pm I}\right)\times\left(\C\setminus\{ 0\}\right)$, 
where 
the action of the element $c$ of the multiplicative group $\C\setminus\{ 0\}$ is 
equal to multiplication by $c$ in the fibers. Via the isomorphisms 
(\ref{Adiso}) and (\ref{LB}), we obtain the equivalent statement that 
{\em the automorphism group of the complex line bundle $Q$ is equal to 
the Cartesian product of the left action of $\op{SO}(3,\,\C )$ and the 
multiplications by nonzero complex numbers on the fibers}. 

Recall that the equivariance of $\Omega$ under the left action of 
$\op{SO}(3,\,\C )$ says that if $\op{L}_g$ denotes the left action of  
$g\in\op{SO}(3,\,\C )$ on 
$Q\setminus 0$, then $\op{L}_g^*\Omega =g\,\Omega$, where in the right 
hand side we let act $g$ on the values in $\C ^3$ of $\Omega$. 
On the other hand the multiplication by $c\in\C\setminus\{ 0\}$ in the 
fibers of $Q\setminus\{ 0\}$ sends $\Omega$ to $c^{-1}\,\Omega$, 
cf. (\ref{Omega}). 
\label{autoQrem}
\end{remark}

\begin{remark}
Let $M$ be a minimal surface in $\R ^3$ and let $C=n'(M)$ be the corresponding 
complex analytic curve in $Q$. The left action of $g\in\op{SO}(3,\,\C )$ maps 
$C$ to the analytic curve $\op{L}_g(C)$, and we obtain a corresponding 
a minimal surface $_gM$ in $\R ^3$ by integrating $\omega$ 
over $\op{L}_g(C)$. If  $g\in\op{SO}(3)$, then it follows from the 
$\op{SO}(3)$\--equivariance of $\omega$ that $_gM=g(M)$, which 
is obtained from $M$ by applying the  
rotation $g$ in $\R ^3$.  (As usual, we work modulo 
translations of the minimal surface.) 
If  $g\in\op{SO}(3,\,\C )$ but $g\notin\op{SO}(3)$, then the 
relation between $M$ and $_gM$ is less straighforward. The best one can 
say is that if  $\Gamma$ is the isotropic complex analytic curve in $\op{C}^3 $ such that 
$M=\op{Re}\Gamma$, where $\Gamma$ is obtained by  means of integration of $\Omega |_C$, 
then the $\op{SO}(3,\,\C )$\--equivariance of  $\Omega$ implies that 
$_gM=\op{Re}(g(\Gamma ))$. 

The action of  multiplication by $c\in\C\setminus\{ 0\}$ on the fibers of 
$Q$ maps the 
complex analytic curve $C$ in $Q$ to the complex analytic curve $c\, C$ in $Q$. 
If $c\in\R$, 
then the minimal surface which is obtained by 
integrating $\omega =q^{-1}\,\op{d}\! s$ over $c\, C$ is equal to the
$c^{-1}\, M$, a homothetic or reflected homothetic image 
of the minimal surface $M$ if $c>0$ or $c<0$, respectively. 

If on the other hand $c=\op{e}^{\fop{i}\phi}$, $\phi\in\R$, 
then the minimal surfaces obtained by integrating $\omega$ over $c\, C$ form a 
circle of minimal surfaces, which is called the {\em family of minimal surfaces 
which is associated to $M$}, cf. DHKW \cite[p. 96, 97]{dhkw}. 
If $\Gamma$ is the isotropic 
complex analytic curve in $\C ^3$ such that $M=\op{Re}\Gamma$, then 
the minimal surfaces associated to $M$ are the 
$\op{Re}\left(\op{e}^{-\fop{i}\phi}\,\Gamma\right)$, $\phi\in\R$. 
The minimal surface $\op{Re}(\op{i}\,\Gamma )$ 
is called the {\em adjoint surface} of $M$, cf. DHKW 
\cite[p. 91]{dhkw}.  
\end{remark}

\begin{remark}
A natural compactification $\overline{Q}$ (in algebraic geometry called completion) of 
$Q$ arises as the ${\C}{\proj}^1$ bundle over $\op{SL}(2,\,\C )/L\simeq{\C}{\proj}^1$,  
which is defined by 
\begin{equation}
\overline{Q}=\op{SL}(2,\,\C )\times _L{\C}{\proj}^1,
\label{Qbar}
\end{equation}
in which $\left(\begin{array}{cc}a&0\\b&1/a\end{array}\right)\in L$ acts on 
${\C}{\proj}^1$ by sending $\C\, (q,\, r)$ to $\C\,\left( a^4\, q,\, r\right)$, 
which is the projective action of the linear transformation  
$\left(\begin{array}{cc}a^4&0\\0&1\end{array}\right)$ on $\C ^2$. 
Actually, the second order contact bundle of minimal surfaces in $\R ^3$,
modulo translations, is 
naturally defined as $\overline{Q}$ rather than as $Q$. The addition 
to $Q$ of the points 
of $\overline{Q}\setminus Q$ at infinity corresponds to
the allowance of singularities of the 
minimal surface at which the Gauss curvature tends to $-\infty$.

The bundle $\overline{Q}$ can be identified with the 
{\em fourth Hirzebruch surface}, cf. Barth, Peters and Van de Ven 
\cite[p. 141]{bpv}.  
In Hirzebruch \cite{h}, the complex surface $\Sigma _n$ 
has been introduced as the algebraic subvariety 
of ${\C}{\proj}^1\times {\C}{\proj}^2$, which is 
defined by the equation $q_1\, {u_1}^n-q_2\, {u_2}^n=0$, 
if $\left[ u_1:\, u_2\right]$ and $\left[ q_0:\, q_1:\, q_2\right]$ 
denote the projective coordinates in ${\C}{\proj}^1$ and ${\C}{\proj}^2$, 
respectively. The projection to $\left[ u_1:\, u_2\right]$ exhibits 
$\Sigma _n$ as a ${\C}{\proj}^1$ bundle over ${\C}{\proj}^1$. 
We may regard the subset where $q_1=q_2=0$ as the section at infinity.  
If we delete this, then we obtain a complex line bundle over 
${\C}{\proj}^1$.  
Over the $\left[ u_1:\, u_2\right]$ with $u_1\neq 0$, 
this complex line bundle has the trivialization 
$(u,\, q)\mapsto \left([1:u],\,\left[ q:\, u^n:\, 1\right]\right)$, 
and over the $\left[ u_1:\, u_2\right]$ with $u_2\neq 0$, 
it has the trivialization 
$(v,\, r)\mapsto \left([v:1],\,\left[ r:\, 1:\, v^n\right]\right)$. 
The image points are the same if and only if 
$v=u^{-1}$ and $r=u^{-n}\, q$, which is the 
retrivialization (\ref{uv}) if $n=4$.  

Because there is a polynomial embedding from 
${\C}{\proj}^1\times {\C}{\proj}^2$ onto  
an algebraic subvariety of ${\C}{\proj}^4$, cf. Shafarevich 
\cite[p. 43]{sha},  
it follows that $\overline{Q}\simeq\Sigma _4$ is a complex projective 
variety. 
\label{Qbarrem}
\end{remark}

\begin{remark}
The identification (\ref{Xx}) of  $\C ^3$ with $\mbox{\gothic sl}(2,\,\C )$ 
has been chosen such that we have $x\in\R ^3$ if and only if $\Xi (x)$ is an 
anti\--selfadjoint $2\times 2$\--matrix. These matrices form the Lie algebra 
of the special unitary group $\op{SU}(2)$, the group of unitary $2\times 2$\--matrices 
with determinant equal to one. The group 
$\op{SU}(2)$ is a so\--called compact real
form of $\op{SL}(2,\,\C )$. We have 
that $\op{SU}(2)$ is diffeomorphic to the unit sphere in $\R ^4$, hence simply 
connected, and the two\--fold covering $\op{Ad}:\op{SU}(2)\to\op{SO}(3)$ is 
the classical way of exhibiting $\op{SU}(2)$ as the universal 
covering of $\op{SO}(3)$, cf. \cite[Sec. 1.2.B]{dk}. 

$L\cap\op{SU}(2)$ is equal to the group $\op{U}(1)$ of the matrices $g=\left(\begin{array}{cc}
a&0\\0&1/a\end{array}\right)$ such that $a\in\C$ and $|a|=1$. If $a=\op{e}^{\fop{i}\psi}$, then $\op{Ad}g$ is 
equal to the rotation about the vertical axis through the angle $\phi =2\psi$. 
It follows that $\op{Ad}\op{U}(1)=\op{SO}(2)$, where the doubling of the angle expresses that 
the adjoint representation is a two\--fold covering. The homomorphisms 
$\op{Ad}:\op{SU}(2)\to\op{SO}(3)$ and $\op{Ad}:\op{U}(1)\to\op{SO}(2)$ 
induce an isomorphism from $\op{SU}(2)/\op{U}(1)$ onto 
$\op{SO}(3)/\op{SO}(2)\simeq S$, 
which leads to still another model for the unit sphere $S$. 
\end{remark} 

\section{Some Applications}
In this section we discuss some applications of the correspondence between 
minimal surfaces $M$ in $\R ^3$ and complex analytic curves $C$ in $Q$. 
We will restrict our attention to the case that $M$ is smoothly immersed in 
$\R ^3$, although it would have been natural to also allow branch points, 
even with infinite curvature 
(in which case $C$ is a curve in the compactification 
$\overline{Q}$ of $Q$). 

\subsection{Flat Points}
\label{fpsubsec}
Let $M$ be a smoothly immersed, connected and  non\--flat minimal surface in $\R ^3$. 
Then the set $M_0$ of points in $M$ where the curvature vanishes is a closed subset of $M$, 
consisting of isolated points. Because the Gauss map $n=\pi\circ nÒ$ 
is a local diffeomorphism on $M\setminus M_0$, it follows that the mapping 
$n'$ is an 
immersion from $M\setminus M_0$ to a smooth (immersed) analytic 
curve $C$ in $Q\setminus 0$ which intersects the fibers of the 
projection $\pi :Q\to S$ transversally. 

If we use the orientation on $M$ which makes 
the Gauss map orientation preserving and provide $M$ with the 
corresponding complex structure, cf. Section \ref{weierstrasssubsec}, then 
$n'$ is a complex analytic mapping from $M$ to $Q$, where 
$x_0\in M_0$ if and only if $n'\left( x_0\right)$ belongs to the 
zero section of $Q$. Let $x_0\in M_0$ and let $z\mapsto x(z)$, 
$z\in D$ be a holomorphic (= isothermal) coordinate in an open neighborhood 
of $x_0$ in $M$. Here $D=\{ z\in\C\mid |z|<\epsilon\}$ is a disk 
around the origin in the complex plane and $x(0)=x_0$, which means that 
$z=0$ corresponds to the point $x_0$. By means of a rotation in $\R ^3$ 
we can arrange that 
$n\left( x_0\right) =e_3$. If we use the $(u,\, q)$\--trivialization of $Q$ over  
$S\setminus\{-e_3\}$ as in (\ref{uOmega}), 
cf. also Remark \ref{stereographrem}, 
then $z\mapsto n'(x(z))$ corresponds to a pair of complex analytic functions 
$z\mapsto u(z)$, $z\mapsto q(z)$ on $D$, where $u(z)$ is the $u$\--coordinate 
of the Gauss map. Note that $u(0)=0$ because $n\left( x_0\right) =e_3$ and $q(0)=0$ 
because the curvature of $M$ is equal to zero at $x_0$. 
The pull\--back of $\Omega$ under the mapping $z\mapsto n'(x(z))$ 
is equal to 
\begin{equation}
\Omega |_C=\frac{-2}{q(z)}\, u'(z)\,\op{d}\! z\, 
\left(1-u(z)^2,\,\op{i}\,\left( 1+u(z)^2\right),\, 2u(z)\right) ,
\label{dxdz}
\end{equation}
cf. (\ref{uOmega}). Recall that $M$ is obtained 
from $C$ by means of integration of $\omega =\op{Re}\Omega$ over $C$, 
hence by means of integration of the real part of (\ref{dxdz}) 
over the disc $D$ in the $z$\--plane. 

\begin{remark}
The formula (\ref{dxdz}) can be identified with the {\em Enneper\--Weierstrass 
representation formula} as in DHKW \cite[(7) on p. 108]{dhkw}, 
if we take $\nu (z)=u(z)$ and $\mu (z)= -4u'(z)/q(z)$. 
However in contrast with our situation, it is assumed in 
the Enneper\--Weierstrass representation formula that 
$\mu$ is holomorphic, $\nu$ is meromorphic, and 
$\mu\,\nu ^2$ is holomorphic on a simply connected domain in $\C$. 
\end{remark}

Because $M$ is not flat, the function $u(z)$ is not equal to a constant and 
also $q$ is not constantly equal to zero. It follows that there are positive 
integers $k$ and $l$ such that 
\begin{equation} 
0\leq i<k\Rightarrow u^{(i)}(0)=0,\; u^{(k)}(0)\neq 0
\quad\mbox{\rm and}\quad
0\leq j<l\Rightarrow q^{(j)}(0)=0,\; q^{(l)}(0)\neq 0.
\label{kl}
\end{equation}
An equivalent condition is that $u(z)/z^k$ and $q(z)/z^l$ extend to an analytic function 
$\widetilde{u}$ and $\widetilde{q}$ on $D$, respectively, such that  
$\widetilde{u}(0)=u^{(k)}(0)/k!\neq 0$ and 
$\widetilde{q}(0)=q^{(l)}(0)/l!\neq 0$.

Because for finite nonzero $\op{d}\! z$ 
the correponding tangent vector $\op{d}\! x$ of $M$ has a finite and nonzero 
limit as $z\to 0$, we obtain from (\ref{dxdz}) that necessarily 
$l=k-1$. Conversely, if $l=k-1$, then the pull\--back of $\Omega$ by means 
of the mapping $z\mapsto n'(x(z))$ is holomorphic and nonzero at $z=0$, 
which implies that (\ref{dxdz}) defines a smooth piece of a minimal surface. 

The number $k$ in (\ref{kl}) is the local mapping degree of the Gauss map 
at the point $x_0$. For $s\in S$ close to $e_3$ there are $k$ points 
$x$ close to $x_0$ such that $n(x)=s$, and the curvature is not equal to zero 
at each of these points $x$.

It also follows from (\ref{kl}) that the local piece 
of the curve $C$ intersects the zero section $0$ of $Q$ at 
$\left( e_3,\, 0\right)$ 
with multiplicity equal to $l$, and the fiber with multiplicity equal to $k$. 
 In this case the intersection multiplicity 
with the fiber is one more than the intersection multiplicity with the zero 
section, which implies that

\subsection{Flat Points at Infinity}
\label{inftysubsec} 
It follows from the Puiseux theory that a subset $C$ in the $(u,\, q)$\--space 
near the origin is equal to the image of a complex analytic map 
$z\mapsto (u(z),\, q(z))$ which satisfies 
(\ref{kl}), if and only if it is equal to the zeroset of an analytic 
function $f(u,\, q)$ such that $f(0,\, 0)=0$ and $f$ has neither $u$ nor $q$ 
as a factor, cf. {\L}ojasiewicz \cite[p. 173]{l}. In other words, if and only if $C$ in a neighborhood $U$ 
of $\left( e_3,\, 0\right)$ in $Q$ is equal to a {\em one\--dimensional 
complex analytic subset of $Q$}, which contains $\left( e_3,\, 0\right)$ but does not 
contain the intersection with $U$ of the zero section or the fiber $Q_{e_3}$. 

Note that, at $\left( e_3,\, 0\right)$, the intersection number of $C$ with the 
fiber, which is equal to the local mapping degree of the 
restriction to $C$ of the 
projection $\pi :Q\to S$, is equal to $k$, whereas the 
intersection number of $C$ with the zero section is equal to $l$. 
$C$ is smooth at $\left( e_3,\, 0\right)$ if $k=1$ or $l=1$. 
We have concluded at the end of Section 
\ref{fpsubsec} that $\left( e_3,\, 0\right)$ corresponds to a finite 
limit point $x_0$ of a minimal surface $M$ if and only if $l\leq k-1$, with 
equality if and only if $M$ is smooth at $x_0$.   
Note that $l=k-1$ implies that  
$C$ is tangent to the fiber, even if $k\geq 3$, when  
$C$ has a singular point at $\left( e_3,\, 0\right)$. 

We assume from now on that $l\geq k$, and study the behaviour  
of the immersed minimal surface $M$ in $\R ^3$ which 
is obtained by integration of 
$\omega =\op{Re}\Omega$ over $C\setminus\{\left( e_3,\, 0\right)\}$.  
If $l\geq k$, then (\ref{dxdz}) has 
a convergent $\C ^3$\--valued Laurent series of the form 
$\sum_{i\geq 0}\,  z^{-l+k-1+i}\,\op{d}\! z\, c_i$, with 
$\C ^3$\--valued coefficients $c_r$ and 
\begin{equation}
c_0=-\frac{2k\,\widetilde{u}(0)}{\widetilde{q}(0)}\, (1,\,\op{i},\, 0)\neq 0. 
\label{c0}
\end{equation}
Its integral therefore has a 
convergent expansion of the form 
\begin{equation}
a+\sum_{i=0}^{l-k-1}\,\frac{z^{-l+k+i}}{-l+k+i}\, c_i
+(\log z)\, c_{l-k}+\sum_{j=l-k+1}^{\infty}\, 
\frac{z^{-l+k+j}}{-l+k+j}\, c_j,
\label{series}
\end{equation}
where $a\in\C ^3$ is a constant, and the finite sum of the poles is absent if $l=k$. 
In particular the immersed minimal surface $M$, obtained by 
integrating the real part of (\ref{dxdz}), tends to infinity 
when $l\geq k$. 

The multi\--valuedness of the logarithm leads to a period vector 
\begin{equation}
p=\op{Re}\left( 2\pi\op{i}\, c_{l-k}\right)\in\R ^3, 
\label{p}
\end{equation}
which is a nonzero vector in the horizontal plane when $l=k$. 
If $l>k$, then it can happen that $p=0$. 

Let $l>k$. If $z$ remains in a sector $0<|z|<\epsilon$,  
$\alpha _1<\op{arg}z<\alpha _2$, then the term in (\ref{series}) with $i=0$ 
dominates, and it follows that the image is asymptotically 
equal to the graph of a function $x_3=f(y)$ 
over a domain in the $y=\left( x_1,\, x_2\right)$\--plane which 
asymptotically (for $\epsilon\downarrow 0$) is close to a 
sector with $\| y\|$ larger than a constant times $(1/\epsilon )^{l-k}$ 
and angle $(l-k)\,\left(\alpha _2-\alpha _1\right)$. If 
$(l-k)\,\left(\alpha _2-\alpha _1\right) >2\pi$, then the 
function $f$ has to be interpreted as a multi\--valued function 
on the overlap. The gradient of $f$ converges to zero 
if  $\| y\|\to\infty$ in the sector. 

If the period vector $p$ in (\ref{p}) is equal to zero, then 
$f(y)$ returns to the same value after $y$ has made $l-k$ turns 
in the sector, which forces $M$ to have self\--intersections along 
real one\--dimensional curves if $l-k>1$. If $p\neq 0$, then a local piece 
$U$ of $M$ will return to $U+p$ when $y$, with large $\| y\|$, 
has made $l-k$ turns in the plane. The number $l-k$ is called 
the {\em spinning number of the end at infinity} of $M$ if 
$p=0$, cf. Hoffman and Karcher 
\cite[(2.16) on p. 20]{hk}. We will use this name for $l-k$  
also in the case that the period vector $p$ is not equal to zero. 
Note that {\em if there is no periodicity, then the end at 
infinity is embedded, 
if and only if the spinning number is equal to one}. 

The leading term in 
the asymptotic expansion for $f(y)$ as $y$ runs to infinity 
in the sector is equal to a nonzero constant times $\log\| y\|$,  
if for each $0<i<l-k$ the third component of $c_r$ is equal to zero, 
but the real part of the third component of $c_{l-k}$ is not equal 
to zero. In all other cases, $f(y)$ converges to a constant 
as $y$ runs to infinity in the sector. In the literature, where one 
assumes that the period vector $p$ is equal to zero,  the end is called 
{\em flat} or {\em planar} if $f(y)$ converges to a constant, 
and {\em of catenoid type} if $f(y)$ has logarithmic 
growth, cf. DHKW \cite[p. 198]{dhkw}. 

If $l=k$, (when the period vector $p$ in (\ref{p}) is always 
nonzero) then the minimal surface image of the sector $
0<|z|<\epsilon$,  $\alpha _1<\op{arg}z<\alpha _2$ 
lies over a subset of the $y$\--plane which is 
asymptotic to a half strip with width equal to 
$|c_0|\,\left(\alpha _2-\alpha _1\right)$, at a distance to the 
origin equal to $|c_0|\,\log (1/\epsilon )$.  
Over it, the function $f(y)$ converges to the real part 
of the third component of $a$ as $y$ runs to infinity. 

\bigskip\noindent
Now assume conversely that $M$ is a smoothly immersed nonflat minimal 
surface in $\R ^3$. For $M$, we introduce the following conditions. 
\begin{itemize}
\item[i)] 
$M$ is given by means of a smooth conformal immersion from the punctured 
disc $D_{\rho}\setminus\{ 0\}$ to $\R ^3$, where 
$D_{\rho}=\{ z\in\C\mid |z|<\rho\}$. 
The orientation of $M$ can be chosen such that, if 
$n(x)$ denotes the corresponding normal to $M$ at $x$, the 
pull\--back $z\mapsto n(\epsilon (z))$ of the Gauss mapping 
by means of $\epsilon$ is a complex analytic mapping from 
$D_{\rho}\setminus\{ 0\}$ to the Riemann sphere $S$. 
It follows that $n'\circ\epsilon$ is a complex analytic mapping from 
$D_{\rho}\setminus\{ 0\}$ to $Q$, cf. Proposition \ref{ccprop}.  
\item[ii)]  
In order to allow for periodicity, as should be done according 
to DHKW \cite[p. 196]{dhkw}, we allow that $\epsilon$ is 
multi\--valued. However, we will assume that if $x=\epsilon (z)$ has 
moved to $\widetilde{x}$ when $z$ has run around the origin once in 
$D_{\rho}\setminus\{ 0\}$, in the positive direction, then 
$n(x)=n(\widetilde{x})$ and $n'(\widetilde{x})=n'(x)$. 
In other words, if we write $p=\widetilde{x}-x$, then 
the translate $M+p$ of $M$ 
osculates $M$ at $\widetilde{x}$. This assumption is equivalent to 
the condition that $z\mapsto n'(\epsilon (z))$ is a single\--valued 
complex analytic mapping from $D_{\rho}\setminus\{ 0\}$ to $Q$. In the latter 
case the integral of $\omega$ over the image  has $p$ as a period, 
which means that actually $M+p=M$.  
\item[iii)] There exists $0<\sigma\leq\tau$ 
such that the image 
$\{ n(\epsilon (z))\mid z\in D_{\sigma}\setminus\{ 0\}\}$ 
of the Gauss map 
misses at least three points of $S$. 
\item[iii')] For any $0<\sigma <\rho$, the integral of 
the Gaussian curvature $K$ over $\epsilon\left( D_{\sigma}\right) /\Z\, p$ 
is finite. 
\item[iv)] If $\gamma :[0,\,\infty [\to D_{\rho}\setminus\{ 0\}$ 
is a $\op{C}^1$ curve such that $\lim_{t\to\infty}\,\gamma (t)=0$, 
then the curve $t\mapsto\epsilon(\gamma (t))$ in $\R ^3$ 
has infinite Euclidean length. 
\end{itemize}

\begin{lemma}
The following conditions {\em a)}, {\em b)}, and {\em c)} are 
equivalent. 
\begin{itemize}
\item[{\em a)}] Let $C_{\sigma}$ denote the image  
$n'\circ\epsilon\left( D_{\sigma}\setminus\{ 0\}\right)$ of 
the minimal surface $M_{\sigma}=\epsilon\left( D_{\sigma}\setminus\{ 0\}\right)$ 
under the mapping $n':M_{\sigma}\to Q$. 
Then there exists $0<\sigma<\rho$ such that $C_{\sigma}$, 
after adding a point $q_0$ on the zero section of $Q$,  
is equal to the germ at $q_0$ of a complex one\--dimensional analytic 
subset $C$ of $Q$, where the intersection number $k$ of $C$ at $q_0$ 
with the fiber is less than or equal to the intersection number 
$l$ of $C$ at $q_0$ with the zero section. 
\item[{\em b)}] {\em i) \& ii) \& iii') \& iv)}. 
\item[{\em c)}] {\em i) \& ii) \& iii) \& iv)}.
\end{itemize} 
\label{flatinftylem}
\end{lemma}
\begin{proof}
For the a) $\Longrightarrow$ b) we only need to prove that (iii') follows 
from the description of $C$. Note that $M$ is obtained 
by means of integration of $\omega$ over $C\setminus\{ q_0\}$, 
which implies that $n(\epsilon (z))$ converges to $q_0$ when $z\to 0$. 
The Gaussian curvature $K$ is equal to minus the Jacobi determinant of the 
Gauss map $n:M\to S$, if we provide $M$ and $S$ with the Euclidean non\--oriented 
area forms $\op{d}\! _2x$ and $\op{d}\! _2s$, respectively. Furthermore 
the restriction to $C$ of the projection $\pi :Q\to S$ has 
local mapping degree equal to $k$. It therefore follows that 
\begin{equation}
\int_{n^{-1}(U)\cap M/\Z\, p}\, K(x)\,\op{d}\! _2x=\, -k\,\op{area}(U), 
\label{intKU}
\end{equation}
if $U$ is a small open neighborhood of $q_0$ in $S$ and we take the pre\--image 
$n^{-1}(U)$ of $U$ under the Gauss map in $M$. If 
the period vector $p$ in (\ref{p}) is nonzero, then $M$ 
is taken modulo the translates over integral multiples of $p$. 
Because $\lim _{z\to 0}\, n(\epsilon (z)) =q_0$, there 
exists $0<\tau <\sigma$ such that 
$D_{\tau }\setminus\{ 0\}(n\circ\epsilon )^{-1}(U)$, 
which in view of (\ref{intKU}) implies that the integral of $K$ over 
$\epsilon\left( D_{\tau}\setminus\{ 0\}\right) /\Z\, p$ is finite. 
Because the integral of $K$ over 
$\epsilon\left( D_{\sigma}\setminus D_{\tau}\right) /\Z\, p$ is 
obviously finite, the conclusion iii') follows. 

If iii') holds then the integral of $K$ 
over $\epsilon\left( D_{\sigma}\setminus\{ 0\}\right) /\Z\, p$ 
converges to zero as $\sigma\to 0$. The argument for 
(\ref{intKU}) implies that the area of 
$n\circ\epsilon\left( D_{\sigma}\setminus\{ 0\}\right)$ 
is less than or equal to $-k$ times the integral of 
$K$ over $\epsilon\left( D_{\sigma}\setminus\{ 0\}\right) /\Z\, p$, 
which implies that there exists $0<\sigma \leq\tau$ such that 
the area of 
$n\circ\epsilon\left( D_{\sigma}\setminus\{ 0\}\right)$ 
is strictly smaller than $4\pi$, the area of $S$. 
This implies iii), and therefore we have proved that b) $\Longrightarrow$ c). 
 
Now assume that c) holds. Using the great Picard theorem, 
cf. Conway \cite[p. 300]{conway}, we obtain from iii) 
that $n(\epsilon (z))$ extends to 
a complex analytic mapping from $D$ to $S$. 
(I learned this argument from Osserman \cite[p. 397]{o63}.) 
By means of a rotation 
in $\R ^3$, we can arrange that the limit point 
$\lim_{z\to 0}\, n(\epsilon (z))$ is equal to $e_3$. 
Using the trivialization of the bundle $Q$ 
over $S\setminus\{ -e_3\}$ by means of the $(u,\, q)$\--coordinates 
in (\ref{uOmega}), the mapping $z\mapsto n'(\epsilon (z))
:D\setminus\{ 0\}\to Q$ is represented by two complex analytic 
functions $u(z)$ and $q(z)$ on $D\setminus\{ 0\}$, 
where $u(z)$ extends to a complex analytic 
function on $D$ such that $u(0)=0$. 

If $0<\sigma <\rho$, then there is a $\delta >0$ such that, for each 
$z\in D_{\sigma}\setminus\{ 0\}$, the distance $d(\epsilon (z))$ 
from $x=\epsilon (z)$ to the boundary of $M$ is larger than 
or equal to $\delta$. Here $d(x)$ is defined as the infinmum of 
the Euclidean length of smooth curves in $M$ which start at 
$x$ and leave every compact subset of $M$. 
A theorem of Osserman \cite{o59} implies that if in a part 
$M_0$ of a smoothly immersed minimal surface the Gauss map 
avoids a neighborhood of a point on $S$, and the distance 
to the boundary of $M$ is bounded in $M_0$, then the 
Gaussian curvature $K$ is bounded on $M_0$. In view of 
(\ref{Kab}), it follows therefore from iii) and iv) 
that the function $q(z)$ is bounded on $D_{\sigma}\setminus\{ 0\}$, 
and therefore it extends to a complex analytic function 
on $D$, which we denote with the same letter. 

The minimal surface immersion $\epsilon (z)$ is found back from 
$(u(z),\, q(z))\in Q$ by 
means of integration of the real part of (\ref{dxdz}). 
Because $\epsilon (z)$ cannot have a finite limit 
as $z\to 0$ in view of iv), it follows that $q(0)=0$ 
and that (\ref{kl}) holds with $l\geq k$. This shows 
that c) $\Longrightarrow$ a) and we have proved the 
equivalence of a), b) and c). 
\end{proof}

We will say that $M$ has a 
{\em flat point at infinity} if any of the equivalent conditions 
a), b), c) in Lemma \ref{flatinftylem} is satisfied. 
(We apologize for 
the collision of this terminology  
with the distinction between flat and catenoid ends at infinity.  
In our terminology both types of ends correspond to flat points 
at infinity.) 
Lemma \ref{flatinftylem} says that the condition iii) can be replaced 
by iii'). Moreover, if $M$ satisfies i), ii), iii) and iv), then  
Lemma \ref{flatinftylem} implies that we have the seemingly much stronger 
asymptotic expansions for $M$ as described after (\ref{series}).

\subsection{Algebraic Curves}
\label{algsubsec}
We now turn to a global version of the conditions i), ii) in 
Section \ref{inftysubsec}. We define a smoothly immersed minimal surface 
in $\R ^3$ {\em of finite topological type} as a multi\--valued 
smooth minimal surface immersion $\epsilon :D\setminus D_{\infty}\to\R ^3$, 
with the following 
properties. 
\begin{itemize}
\item[i)] $D$ is a compact Riemann surface and $D_{\infty}$ is a finite 
subset of $D$, possibly empty.  
If we provide the immersed minimal surface $M$ with the complex structure 
such that the Gauss mapping is a complex analytic mapping $n$ from $M$ 
to the Riemann sphere $S$, then $n\circ\epsilon$ is a complex analytic mapping 
from $D\setminus D_{\infty}$ to $S$. 
In turn this implies that $n'\circ\epsilon$ is a complex analytic 
mapping from $D\setminus D_{\infty}$ to $Q$, cf. Proposition \ref{ccprop}. 
\item[ii)] In order to allow for periodicity, the mapping $\epsilon$ may 
be multi\--valued in the following way. If $\gamma :[0,\, 1]\to D\setminus 
D_{\infty}$ is a curve in $D\setminus D_{\infty}$ such that 
$\gamma (1)=\gamma (0)$ and $\epsilon (\gamma (t))$ has moved from 
$x$ to $\widetilde{x}$ as $t$ runs from $0$ to $1$, then 
$n(\widetilde{x})=n(x)$ and $n'(\widetilde{x})=n'(x)$. This assumption 
is equivalent to the condition that $n'\circ\epsilon$ is a single\--valued 
mapping from $D\setminus D_{\infty}$ to $Q$. The image 
$n'\circ\epsilon\left(D\setminus D_{\infty}\right)$ therefore is 
a complex analytic curve $C$ in $Q$ and the minimal surface 
$M$ is obtained by means of integration of 
$\omega |_C=\op{Re}\Omega |_C$. If 
${\cal P}$ denotes the corresponding set of periods, which is 
an additive subgroup of $\R ^3$, and $\phi =\int _C\omega$ is 
the reconstruction map $:C\to M/{\cal P}$ of Section \ref{recsubsec}, 
then $\phi\circ (n'\circ\epsilon )$ is a single\--valued mapping 
from $D\setminus D_{\infty}$ to $\R ^3/{\cal P}$, which is the 
single\--valued realization of $\epsilon$. Note that ${\cal P}$ 
is generated by the $\langle\left[\gamma _i\right],
\,\left[ (n'\circ\epsilon )^*\omega\right]\rangle$, where 
$\left[\gamma _i\right]$ runs over a set of generators of 
the image of $\op{H}_1(D\setminus D_{\infty},\,\Z )$ in 
$\op{H}_1(D\setminus D_{\infty},\,\R )$. 
\end{itemize}

If $C$ is a germ at $q$ of an analytic subset of $Q$, then a
{\em component of $C$ at $q$} is defined as a  
set of the form $A\cup\{ q\}$, in which $A$ is a connected component of 
$B\cap (C\setminus\{ q\})$ and $B$ is a small ball around $q$ in $Q$. 
We also recall the compactification $\overline{Q}$ of $Q$, 
the ${\C}{\proj}^1$\--bundle $\Sigma _4$ over ${\C}{\proj}^1$, 
introduced in Remark \ref{Qbarrem}.

\begin{proposition}
Let $\epsilon :D\setminus D_{\infty}\to M/{\cal P}$ be a smoothly immersed 
nonflat minimal surface in $\R ^3$, modulo its periods, of finite topological 
type. Assume moreover that, for every $e\in D_{\infty}$ and small disc 
$D(e)$ around $e$ in $D$, $\epsilon |_{D(e)}$ is a flat point at infinity. 
Then, if we add to $n'\circ\epsilon\left( D\setminus D_{\infty}\right)$ 
the finitely many limit points $q(e)$, $e\in D_{\infty}$, on the zero 
section $0=0_{_{Q}}$ of $Q$, we obtain a complex {\em algebraic} curve $C$ in $\overline{Q}$ 
with the following properties.
\begin{itemize}
\item[{\em a)}] $C\subset Q$. 
\item[{\em b)}] For every $q\in C\setminus 0_{_{Q}}$, every component of $C$ at $q$ 
is smooth and transversal to the fiber through $q$ of $\pi :Q\to S$. 
\item[{\em c)}] If $q\in C\cap 0_{_{Q}}$, then each component of $C$ at $q$ 
has intersection number $k$ with the fiber through $q$ and intersection 
number $l$ with $0_{_{Q}}$, where $l\geq k-1$. 
\end{itemize}
The components of $C$ at points of $C\cap 0_{_Q}$ such that $l\geq k$ 
correspond bijectively to the points of $D_{\infty}$, or the flat points at inifinity. 
The others correspond to the finite flat points of $M/{\cal P}$. 

Conversely, if $C$ is an algebraic curve in $\overline{Q}$ such that 
{\em a), b), c)} hold, then it arises from a minimal surface $M/{\cal P}$ 
of finite topological type and with flat points at infinity as above, 
where $M/{\cal P}$ is obtained from $C$ by means of integration of 
$\omega |_C=\op{Re}\Omega |_C$. 
\label{algprop}
\end{proposition}

\begin{proof}
The assumptions on $\epsilon$ imply that $n'\circ\epsilon$ extends 
to a complex analytic mapping from $D$ to $Q$. Its image $C$ is a 
compact subset of $Q$, and it follows from Lemma \ref{flatinftylem} 
that $C$ is an analytic subset of $Q$. 
Because $\overline{Q}$ is a complex projective variety, 
cf. Remark \ref{Qbarrem}, 
and Chow's theorem states that every complex analytic 
subset of a complex projective space is algebraic, 
cf. Griffiths and Harris 
\cite[p. 167]{gh}, the conclusion is that $C$ is a complex 
algebraic curve in $\overline{Q}$, which moreover satisfies a). 
The properties b) and c) follow from the discussion in Section 
\ref{fpsubsec} and \ref{inftysubsec}. 

For the converse, let $\eta :D\to C$ be the normalization of $C$, 
as defined in for instance {\L}ojasiewicz \cite[p. 343, 344]{l}. 
It is obtained by replacing, for every nonsmooth point $q$ of $C$ 
and small ball $B$ in $Q$ around $q$, the set 
$C\cap B$ by the disjoint union of the 
components of $C$ at $q$, where each component is parametrized 
by means of a complex analytic map $z\mapsto (u(z),\, q(z))$, with 
$z$ in a small disc in the complex plane, as in (\ref{kl}). 
Then $D$ is a compact Riemann surface and $\eta :D\to C$ is a 
complex analytic mapping such that, for each $q\in C$, the number of 
elements in the pre\--image $\eta ^{-1}(\{ q\})$ is equal to the 
number of components at $q$ of $C$. Define $D_{\infty}$ as the 
finite set of the points $e\in D$ such that $q=\eta (e)\in C\cap 0_{_{Q}}$ 
and, for a small disc $D(e)$ around $e$ in $D$, $\eta (D(e))$ is equal to 
a component at $q$ of $C$ as in c), with $l\geq k$. 

If $\phi =\int_C\omega$ denotes the integration map from $C\setminus 
\eta\left( D_{\infty}\right)$ to the minimal surface $M/{\cal P}$, 
then $\epsilon =\phi\circ\eta$ is a minimal surface of finite 
topological type with flat points at infinity, and $n'\circ\epsilon$ 
maps $D\setminus D_{\infty}$ to $C$ because $n'$ is equal to the 
inverse of $\phi$.  
\end{proof}

\begin{remark}
It is a theorem of Osserman \cite[p.81, 82]{o}, 
%\cite[***]{lawson}
that a complete immersed minimal surface $M$ with finite total curvature 
is of finite topological type. The condition that $M$ itself 
(not some quotient by non\--trivial translations) has finite 
total curvature implies that ${\cal P}=\{ 0\}$. It follows from 
Proposition \ref{algprop} that  
these $M$ correspond to the complex algebraic curves $C$ in 
$\overline{Q}$ which satisfy a), b), c), and moreover have the 
property that all the integrals of $\omega =\op{Re}\Omega$ 
over closed curves in $C\setminus 0_{_{Q}}$ are equal to zero. 
The latter is a very severe restriction on the algebraic curves $C$. 
It remains to be seen in how far a systematic search in the 
parameter range of the complex algebraic curves can be performed in 
order to find (or classify) the examples in which one is interested.   
\label{noperiodrem} 
\end{remark}

\begin{remark}
If one of the conditions a), b), c) for the complex algebraic 
curve $C$ is not satisfied, then the corresponding minimal surface 
is no longer smoothly immersed and has singularities of 
branch point type. If $C$ intersects the section 
$\infty _{_{Q}}:=\overline{Q}\setminus Q$ at infinity, then 
the curvature tends to infinity if one approches the 
singular point, in all other cases the curvature has 
a finite limit at the branch point. 
\label{singrem}
\end{remark}

Let $(u,\, q)$ be the trivialization of $Q\setminus Q_{-e_3}$ over 
$S\setminus\{ -e_3\}$ as in (\ref{uOmega}). Because of the polynomial 
embedding of $\overline{Q}$ in a complex projective space as 
in Remark \ref{Qbarrem}, any complex algebraic curve $C$ in 
$\bar{Q}$ is equal, within $Q\setminus Q_{-e_3}$, to the zeroset of 
a polynomial $f(u,\, q)$ in the two complex variables $u,\, q$, 
which we can write in the form 
\begin{equation}
f(u,\, q)=\sum_{j=0}^d\, f_j(u)\, q^j.
\label{Pj}
\end{equation}
Here, for each $0\leq j\leq d$, $f_j$ is a polynomial 
in one variable and $f_d(u)$ 
is not identically equal to zero. 

The condition that $C\subset Q$ implies that $f_d(u)$ has no 
zeros, which means that it is equal to a nonzero constant, which 
we can take equal to $-1$, because multiplication of $f$ with 
a nonzero constant does not change its zeroset. However, 
the condition $C\subset Q$ implies that the same is true 
for the polynomial equation which is obtained after the 
retrivialization (\ref{uv}). We have 
\begin{equation}
f\left( v^{-1},\, v^{-4}\, r\right) 
=\sum_{j=0}^d\, f_j\left( v^{-1}\right)\, v^{-4j}\, r^j
=v^{-m}\,\sum_{j=0}^d\, v^{m-4j}\, f_j\left( v^{-1}\right)\, r^j,  
\label{fv}
\end{equation}
in which $m$ denotes the maximum of the numbers $4j+\op{degree}f_j$, 
$0\leq j\leq d$. Because it is required that $m-4d=0$, it follows 
that the closure in $\overline{Q}$ of the zeroset of $f$ is contained in $Q$, 
if and only if (\ref{Pj}) holds with 
\begin{equation}
f_d=\, -1,\quad \op{deg}f_j\leq 4(d-j)\quad\mbox{\rm for}\quad 0\leq j\leq d-1.
\label{Pjcond}
\end{equation} 

We have the identities $d$ = the intersection number of $C$ with each fiber 
of $\pi :Q\to S$ = the degree of the mapping 
$\pi |_C:C\to S$ = the degree of the 
Gauss map from $M/{\cal P}$ to $S$. (Strictly speaking, here $M/{\cal P}$ 
should be replaced by its compactification, obtained by the addition of a 
finite set of points at infinity. Note that 
$d$ is also equal to the degree of the 
mapping $n\circ\epsilon$ from $D$ to $S$.) 
In view of (\ref{intKU}), we obtain that the 
total curvature of the minimal surface $M/{\cal P}$ modulo ${\cal P}$ is equal to 
\begin{equation}
\int_{M/{\cal P}}\, K(x)\,\op{d}\! _2x=\, -4\pi\, d, 
\label{intKM}
\end{equation}
because the area of $S$ is equal to $4\pi$. As observed in Lemma \ref{flatinftylem}, 
the condition that the total curvature of $M/{\cal P}$ is finite is 
equivalent to the condition of having only flat points at infinity. 
Therefore, {\em if the 
complete minimal surface modulo ${\cal P}$ is of finite topological type  
and has finite total curvature, then its total curvature is 
equal to $-4\pi\, d$, where $d$ is equal to the degree of the Gauss map from 
$M/{\cal P}$ to $S$}. 

The intersection points of $C\setminus Q_{-e_3}$ with the zero section $0_{_{Q}}$ of $Q$ correspond to the 
zeros $u$ of the polynomial $f_0$ of degree less than or equal to $4d$, whereas (\ref{fv}) 
shows that at $\left( -e_3,\, 0\right)$ $C$ has an intersection with 
$0_{_{Q}}$ of multiplicity equal to  
$4d-\op{degree}f_0$. It follows that {\em if the degree of the 
Gauss map is equal to $d$, then the intersection number of $C$ with 
the zero section of $Q$ is equal to $4d$}. In other words, {\em the number 
of flat points of $M/{\cal P}$ plus its number of  flat points at infinity, each counted 
with multiplicity, is equal to four times the degree of the Gauss map 
$n:M/{\cal P}\to S$}. 
\begin{remark}
A topological explanation for the last statement 
can be given as follows. The Chern number of the complex line bundle 
$Q$ over $S\simeq {\C}{\proj}^1$ is equal to 4, 
which implies that the self\--intersection number $0_{_{Q}}\cdot 0_{_{Q}}$ of 
$0_{_Q}$ is equal to $4$. 

The self\--intersection number 
$F\cdot F$ of a fiber $F$ of $\pi :\overline{Q}\to S$ is equal to zero, because 
two distinct fibers are disjoint. The argument of Griffiths and Harris 
\cite[p. 518]{gh} yields that every 
real two\--dimensional cycle $[C]\in\op{H}_2(\overline{Q},\,\Z )$ can be written as 
$c\, F+d\, 0_{_{Q}}$, for some integers $c$ and $d$. We then obtain from 
\[
C\cdot F=c\, F\cdot F+d\,\, 0_{_{Q}}\cdot F=d
\]
that $C\cdot F =d$ = the degree of the mapping $\pi |_C:C\to S$. 
On the other hand 
\[
C\cdot 0_{_{Q}}=c\, F\cdot 0_{_{Q}} +d\,\, 0_{_{Q}}\cdot 0_{_{Q}} 
=c+4d.
\]

Let $\infty _{_{Q}}=\overline{Q}\setminus Q$ 
denote the section at infinity of $\pi :\overline{Q}\to S$, 
which has been added to $Q$ in order to compactify $Q$. 
Then $0_{_{Q}}\cdot\,\infty _{_{Q}} =0$ and $F\cdot\,\infty _{_{Q}} =1$, and 
therefore 
\[
C\cdot\,\infty _{_{Q}} =c\, F\cdot\,\infty _{_{Q}} 
+d\,\, 0_{_{Q}}\cdot\,\infty _{_{Q}} =c, 
\]
and we arrive at the conclusion that 
\begin{equation}
C\cdot\, 0_{_{Q}} =C\cdot\,\infty _{_{Q}} +4\, C\cdot F
\label{C0}
\end{equation}
for every real two\--dimensional cycle $C$ in $\overline{Q}$. In particular, 
the intersection number of $C$ with the zero section is equal to four times the 
degree of the mapping $\pi |_C:C\to S$, if and only if the intersection number 
of $C$ with the section at infinity is equal to zero. 

Applying  
(\ref{C0}) to $C=\infty _{_{Q}}$, it also follows that 
$\infty _{_{Q}}\cdot\,\infty _{_{Q}} =\, -4$. The negativity 
of this self\--intersection number implies that $\infty _{_{Q}}$ is rigid, 
in the sense that $\infty _{_{Q}}$ can not be deformed to a nearby complex algebraic curve, because 
the intersection number between complex algebraic curves is always $\geq 0$.  
\label{toprem}
\end{remark}

Let $C$ be the algebraic curve in $\overline{Q}$ which corresponds to the 
minimal surface $\epsilon :D\setminus D_{\infty}\to M/{\cal P}$ of finite topological 
type and flat points at infinity, as in Proposition \ref{algprop}. 
Let $\eta :D\to S$ be the complex analytic extension to $D$ of the Gauss map  
$n\circ\epsilon :D\setminus D_{\infty}\to S$. Then $\eta$ exhibits the compact Riemann 
surface $D$ as a branched covering of $S$, where the branch points in $S\simeq 0_{_{Q}}$ 
correspond to the points $q\in C\cap 0_{_{Q}}$. For each component $A$ at $q$ of $C$, 
we have a corresponding branch point $b\in D$, where the $b$ are all distinct. 
We write $k(b)$ and $l(b)$ for the intersection number of $A$ with the fiber through $q$ and 
the zero section, respectively. 

Let $B$ denote the set of all the branch points $b\in D$, note that $B$ is a finite 
subset of $D$ and that $D_{\infty}$ is equal to the set of $b\in B$ such that 
$l(b)\geq k(b)$. Also recall that $\op{spin}(b)=l(b)-k(b)$ is the 
spinning number of the flat point at infinity, as introduced in the second paragraph 
after (\ref{p}). Let $d$ denote the degree of the Gauss map $n\circ\epsilon :D\to S$, 
$r:=\#\left( D_{\infty}\right)$ the number of flat points at infinity and 
$s:=\sum_{b\in D_{\infty}}\,\op{spin}(b)$ the {\em total spin at infinity} of 
the minimal surface. We then have the following formula for the genus $g$ of $D$:    
\begin{equation}
g=1+d-r/2-s/2. 
\label{genusspin}
\end{equation}
\begin{proof}
Because the genus of $S$ is equal to zero, the Riemann\--Hurwitz formula, 
cf. Farkas and Kra \cite[p. 21]{fk}, yields that 
\begin{equation}
g=1-d+\sum_{b\in B}\,\frac{k(b)-1}{2}.
\label{RH}
\end{equation}
Now we split the sum in the right hand side into a sum over the $b\in D\setminus 
D_{\infty}$ and the $b\in D_{\infty}$. If $b\in D\setminus D_{\infty}$, 
which corresponds to a finite flat point, then $k(b)-1=l(b)$, cf. 
the end of Section \ref{fpsubsec}. On the other hand, if 
$b\in D_{\infty}$ then $k(b)-1=l(b)-\op{spin}(b)-1$, because of the definition 
of $\op{spin}(b)$. The formula (\ref{genusspin}) now follows from 
(\ref{RH}), because $\sum_{b\in B}\, l(b)=C\cdot 0_{_{Q}}=4d$. 
\end{proof}

The formula (\ref{genusspin}) agrees with the Jorge\--Meeks formula, 
cf. Jorge and Meeks \cite[p. 210]{jm}, Hoffman and Karcher 
\cite[(2.16) on p. 20]{hk}, 
where one has made the  
assumption that $M$ itself has finite total curvature, 
not some quotient by a nonzero group of translations,  
cf. Remark \ref{noperiodrem}. In \cite{jm} the Jorge\--Meeks 
formula has been proved with the help of the Gauss\--Bonnet 
formula for the total curvature, whereas our argument is 
based on the topology of the bundle $Q$ which encodes the 
second order contact modulo translations.  

If ${\cal P}=\{ 0\}$, then $\chi (M)=2-2g-r$ is the 
{\em Euler characteristic} of $M$, where $g$ is the genus 
of $D$ and $r=\#\left( D_{\infty}\right)$. Combination of (\ref{intKM}) 
and (\ref{genusspin}) therefore leads to the formula 
\begin{equation}
\int_M\, K(x)\,\op{d}\! _2x=2\pi\, (\chi (M)-s)
\label{Kchi}
\end{equation}
for the total curvature, where $s\geq r$ and $s=r$ if and 
only if all ends at infinity are embedded, cf. Jorge and Meeks 
\cite[p. 210]{jm}. 
The inequality that the total curvature is $\leq 2\pi (\chi (M)-r)$ 
is due to Gackstatter \cite{gack}. 

\subsection{Sections}
\label{secsubsec} 
The simplest complex algebraic curves in $Q$, apart from the fibers, 
are the sections, i.e. the curves $C$ which have intersection number 
with the fibers equal to one. This means that $d=1$ in (\ref{Pj}), (\ref{Pjcond}), 
or that in the $(u,\, q)$\--trivialization $C$ is equal to the graph 
$q=f_0(u)$ of a polynomial $f_0$ of degree $\leq 4$. It follows that 
the holomorphic sections of $Q$ form a $5$\--dimensional complex vector 
space. Of the 10 real parameters which determine the corresponding 
minimal surface, 6 are ``essential'', because we may substract the 3 real 
dimensions of $\op{SO}(3)$ and the parameter of the dilations. 

\begin{remark} The space of holomorphic sections of $Q$ 
can be recognized as the representation space for the irreducible 
representation of $\op{SL}(2,\,\C )$ with highest weight equal to four times the 
fundamental highest weight of $\op{SL}(2,\,\C )$, where $Q$ is the complex 
line bundle over $\op{SL}(2,\,\C )/L\simeq {\C}{\proj}^1$ 
in the Borel\--Weil picture, cf. \cite[(4.12.5), (4.12.7)]{dk}.  
\label{BWrem}
\end{remark}

Because the degree $d$ of the projection $\pi _C:C\to S$ is equal to one, we have for every 
$s\in S\setminus\pi\left( C\cap 0_{_{Q}}\right)$ a unique $x\in M/{\cal P}$ such that 
$n(x)=s$, whereas the is no such $x$ if $s\in\pi\left( C\cap 0_{_{Q}}\right)$. 
Also note that $d=1$ is equivalent to the condition that the total curvature 
of $M/{\cal P}$ is equal to $-4\pi$. 

The topological intersection number of $C$ with the zero section is equal to $4$, 
where the intersection points of $C$ with $0_{_{Q}}$ correspond to the zeros 
of $f_0$, counted with multiplicity, and an intersection at $\left( -e_3,\, 0\right)$ 
with multiplicity $4-\op{deg}f_0$, if $\op{deg}f_0<4$. If $q$ is an intersection 
point of $C$ and $0_{_{Q}}$, then the numbers $k$ and $l$ in (\ref{kl}) satisfy $k=1$, 
whereas $l$ is equal to the multiplicity of the intersection (= the order of the 
zero of $f_0$ if $\pi (q)\neq -e_3$). In particular, because $l\geq k$, 
each $q\in C\cap 0_{_{Q}}$ corresponds to a single flat point at infinity, 
where $\pi (q)$ is equal to the limiting 
position of the normal. Also note that $l\geq k$ for every 
$q\in C\cap 0_{_{Q}}$ implies that {\em the minimal surface has no 
finite flat points}. 

Because $C$ is diffeomorphic to the sphere, the sum of homology classes of 
the positively oriented small loops around the points of $C\cap 0_{_{Q}}$ is equal to zero, which 
implies that the sum of the corresponding period vectors, the integrals of 
$\omega$ over these loops, is also equal to zero. It follows that 
we have no periodicity, i.e. ${\cal P}=\{ 0\}$, if we have only one 
intersection point $q$ of $C$ with $0_{_{Q}}$, which then necessarily has multiplicity 
equal to 4. By means of a rotation in $\R ^3$, we can arrange that $\pi (q)=\, -e_3$, 
which implies that $f_0(u)$ is a nonzero constant. The corresponding minimal 
surface in $\R ^3$ is {\em Enneper's surface}, cf. DHKW 
\cite[pp. 144-149]{dhkw}.    

We have already observed in (\ref{c0}) and (\ref{p}) that if $l=k(=1)$, then 
the integral of $\omega$ over a small loop in $C$ around $q$ is equal to a 
nonzero period vector $p$. If $\pi (q)\neq -e_3$ and $q$ has the coordinates $(a, 0)$ in the 
$(u,\, q)$\--trivialization of $Q\setminus Q_{-e_3}$, then 
$f_0(u)=(u-a)\, g(u)$, where $g$ is a polynomial with $\op{deg}g=\op{deg}f_0-1$ and 
$g(a)\neq 0$. Taking $u(z)=u$, $q(z)=f_0(u)$ in (\ref{dxdz}), we arrive at the formula 
\begin{equation}
p=4\pi\,\op{Im}\left(\frac{1}{g(a)}\,\left( 1-a^2,\,\op{i}\left( 1+a^2\right),\, 2a\right)\right) 
\label{pg}
\end{equation}
for the period $p$. 

Therefore the only remaining case without periodicity is that we have two intersection 
points $q_1$ and $q_2$ of $C$ with $0_{_{Q}}$, each of which of multiplicity two. 
By means of a rotation in $\R ^3$ we can arrange that $\pi\left( q_2\right) =\, -e_3$, 
which implies that $f_0(u)=c\, (u-a)^2$ for some $u\in\C$ and $c\in\C\setminus\{ 0\}$. 
A straighforward calculation shows that in this case the period vector $p$ is equal to 
$4\pi$ times the imaginary part of $c^{-1}\,\left( -2a,\, 2\op{i}a,\, 2\right)$, 
(The period vector of $q_2$ is equal to $-p$, because the sum of the 
period vectors is equal to zero.) We have that $p=0$ 
if and only if $a=0$ and $c\in\R$. The corresponding 
minimal surface in $\R ^3$ is the {\em catenoid}, 
cf. DHKW \cite[pp. 135-138]{dhkw}. The {\em helicoid} is the 
adjoint surface of the catenoid, for which $f_0(u)=\op{i}u^2$, 
which belongs to the {\em associated family of the catenoid}, 
for which $f_0(u)=\op{e}^{\fop{i}\theta}\, u^2$,  
cf. DHKW \cite[pp. 138-140]{dhkw}. The period vector is equal to  
$p=-4\pi\, (\sin\theta )\, e_3$ for the members of this family. 

Because the sum of the period vectors is equal to zero, we obtain at most three 
linearly independent ones, and this actually occurs for the generic section $C$ 
(which necessarily has four simple intersection points with $0_{_{Q}}$). In this case 
${\cal P}$ is a discrete additive subgroup of $\R ^3$, $\R ^3/{\cal P}$ is 
diffeomorphic to the three\--dimensional torus, in which $M/{\cal P}$ is smoothly 
immersed. A case at hand is {\em Scherk's surface}, for which 
$f_0(u)=u^4-1$, cf. DHKW \cite[p. 153]{dhkw}. 

On the other hand, it can also happen that the subgroup ${\cal P}$ of $\R ^3$ is not 
discrete. For instance, one can have nonzero $p_1$ and $p_2$ in ${\cal P}$ such that 
$p_2=\alpha\, p_1$, where $\alpha$ is an irrational real number. In such a case 
the minimal surface will be dense in $\R ^3$, which makes it rather senseless to 
make a picture of the full minimal surface. 

{\em Bour's surfaces} are given by $q=c\, u^m$, where $m\in\R$ and 
$c\in\C\setminus\{ 0\}$, cf. DHKW \cite[p. 149]{dhkw}. These correspond 
to minimal surfaces as described in Proposition \ref{algprop}, 
if and only if $m\in\Q$ and $0\leq m\leq 4$. 
In view of the retrivialization (\ref{uv}), the  
replacement of $m$ by $4-m$ corresponds to a switch in the role 
of $e_3$ and $-e_3$, which means that we may restrict our 
attention to the cases that $0\leq m\leq 2$. For $m=0$ 
we have Enneper's surface, and for $m=2$ the catenoid family. 
If $m=a/b$ with $a,\, b\in\Z$, $b>0$ and $\op{gcd}(a,\, b)=1$, 
then the complex algebraic curve $C$ is defined by the equation 
$q^b=c\, u^a$. The degree of the Gauss mapping is equal to $b$, 
but $M/{\cal P}$ is isomorphic to the Riemann sphere with 
two points deleted. (According to (\ref{RH}), 
the genus of $D$ is equal to zero.)  In this sense, although 
$C$ is not a section, Bour's surface is akin to the minimal 
surfaces defined by a section. The period $p$ can be computed 
from (\ref{dxdz}) with $u(z)=z^b$, $q(z)=\widetilde{c}\, zÔ$, 
where $\widetilde{c}^b=c$.  It turns out that 
{\em the period vector is equal to zero, 
unless $a=2$, $b=1$, or $a=b=1$}.  In both of these 
exceptional cases,  
$C$ is a section with two zeros at opposite points of $S$. 
The first case is the catenoid family. In the second case, where 
one zero is simple and the other has multiplicity 3, we have 
for every nonzero value of $c$ that   
$p$ is a nonzero vector in the horizontal plane.  

\subsection{Hyperelliptic Curves}
\label{hypsubsec}
It follows from (\ref{genusspin}) that a minimal surface as in Proposition \ref{algprop} has {\em no flat points 
at infinity}, i.e. $D_{\infty}=\emptyset$, if and only if 
the genus $g$ of $C$ is equal to $d+1$, if $d$ denotes 
the degree of the Gauss map (= the degree of $\pi |_C:C\to S$).  
\begin{remark}  
In Meeks \cite[Thm 7.1]{meeks} the formula $d=g-1$ was obtained 
as a consequence of the Gauss\--Bonnet formula 
\[
-4\pi\, d=\int_{M/{\cal P}}\, K(x)\,\op{d}\! _2x=2\pi\,\chi (M/{\cal P})
=2\pi\, (2-2g).
\]
He also observed that $g=2$ cannot occur: if $g=2$ then 
$d=1$, hence the Gauss map is an isomorphism from $M/{\cal P}$ 
onto $S$, which implies that $g=0$. 
\end{remark}
More explicitly, we have no flat points at infinity if and only if, 
for each $q\in C\cap 0_{_{Q}}$ and each component $A$ at $q$ of $C$, 
we have that 
$l=k-1$, if $k$ is the intersection number of $A$ with the fiber 
and $l$ is the intersection number of $A$ with the zero section $0_{_{Q}}$. 
This implies that $d\geq k\geq 2$. If we take $d=2$, 
which implies that $g=3$, then it follows that 
$C$ is equal to the zeroset of a polynomial $f$ as in (\ref{Pj}), (\ref{Pjcond}), 
with $d=2$, $f_1(u)\equiv 0$, and $f_0$ a polynomial of degree 7 or 8 with 
only simple zeros. (The degree is equal to 7 if and only if one of the points 
of $C\cap 0_{_Q}$ lies over $-e_3$, which always can be arranged by 
means of a rotation in $\R ^3$.) In this case $C$ is equal to the 
{\em hyperelliptic curve of genus 3} defined by the equation 
\begin{equation}
q^2=f_0(u),
\label{hyperell}
\end{equation}
in which the degree of $f_0$ is equal to 7 or 8 and $f_0$ has only simple zeros. 
$C$ is a smooth complex algebraic curve in $\overline{Q}$ and 
the projection $\pi |_C:C\to S$ is the usual branched covering 
of $S\simeq \C\proj ^1$ by the hyperelliptic curve $C$. 
It is a cyclic covering as described in Barth, Peters and 
Van de Ven \cite[p. 42]{bpv}. 

The minimal surface $M/{\cal P}$ modulo the periods 
has eight finite flat points, corresponding to the eight intersection points 
of $C$ with the zero section $0_{_Q}$ of $Q$. 
Via the coefficients of the polynomial $f_0$ in (\ref{hyperell}), 
the set ${\cal C}$ of these hyperelliptic curves $C$ can be identified with 
a nonvoid open subset of $\C ^9$. 

Let $\oint\Omega :[\gamma]\mapsto \oint_{\gamma}\Omega$ be the mapping 
from $\op{H}_1(C,\,\Z )$ to $\C ^3$, which is obtained 
by means of integration of  the $\C ^3$\--valued one\--form 
$\Omega$ over closed loops in $C$. The image ${\cal P}\left(\Omega |_C\right)$ 
of  $\op{H}_1(C,\,\Z )$ in $\C ^3$ is an additive subgroup of $\C ^3$ 
with $6$ generators, because $\op{H}_1(C,\,\Z )\simeq\Z ^6$, 
cf. Farkas and Kra \cite[I.2.5]{fk}. 
If $\op{Re}$ denotes the real linear mapping from 
$\C ^3$ onto $\R ^3$ which assigns to $z\in\C ^3$ its real part, 
then $\op{Re}{\cal P}\left(\Omega |_C\right) ={\cal P}$, the 
group of the integrals of $\omega =\op{Re}\Omega$ over closed loops 
$\gamma$ in $C$. 
Note that the fact that $\Omega |_C$ is holomorphic implies that the 
restriction of $\omega$ to $C$ is smooth (real analytic). 

\begin{lemma}
${\cal P}\left(\Omega |_C\right)$ is a lattice in $\C ^3$, which means that it has a 
$\Z$\--basis which is also an $\R$\--basis of $\C ^3\simeq\R ^6$. 
This implies that $\C ^3/{\cal P}\left(\Omega |_C\right)$ is compact, a real 
6\--dimensional torus. 

The mapping $\int\Omega : q\mapsto\int_{q_0}^q\,\Omega$, where the integration 
is over a curve in $C$ which runs from $q_0$ to $q$, defines a 
complex analytic embedding from $C$ into $\C ^3/{\cal P}\left(\Omega |_C\right)$. 
The image $(\int\Omega )(C)$ is a compact complex analytic curve in 
$\C ^3/{\cal P}\left(\Omega |_C\right)$ and $\int\Omega$ is a complex analytic diffeomorphism  
from $C$ onto its image. 
\end{lemma}
\begin{proof}
Let ${\cal H}^1(C)$ denote the space of holomorphic $(1,\, 0)$\--forms on 
$C$, also called the abelian differentials on $C$. For any compact Riemann 
surface, the complex dimension of ${\cal H}^1(C)$ is equal to the genus $g$ 
of $C$, cf. Farkas and Kra \cite[p. 62]{fk}, which in our case is equal to 3. 
For each closed curve $\gamma$ in $C$ we have the complex linear 
form $\oint_{\gamma}:\theta\mapsto\oint_{\gamma}\theta$ on ${\cal H}^1(C)$. 
The mapping $\oint :[\gamma ]\to\oint_{\gamma}$ is injective from 
$\op{H}_1(C,\,\Z )$ to the dual space  ${\cal H}^1(C)^*$ of ${\cal H}^1(C)$, 
and its image $\Lambda$ is a lattice in ${\cal H}^1(C)^*$, 
cf. Farkas and Kra \cite[III.2.8]{fk}. $\Lambda$ is 
called the {\em period lattice of $C$}. 
This implies that $\op{Jac}(C):={\cal H}^1(C)^*/\Lambda$ 
is a real 6\--dimensional torus, called the {\em Jacobi variety of $C$}. 

The components of  $\Omega |_C$ are equal to 
\[
\Omega _1=-2\left( 1-u ^2\right)\,\op{d}\! u/q, 
\quad\Omega _2=-2\op{i}\left( 1+u^2\right)\,\op{d}\! u/q,\quad 
\mbox{\rm and}\quad 
\Omega _3=-4u\,\op{d}\! u/q, 
\]
cf. (\ref{uOmega}). Because the one\--forms $u^j\,\op{d}\! u/q$ with  
$j=0,\, 1,\, 2$ form a basis of  ${\cal H}^1(C)$, 
cf. Farkas and Kra \cite[p. 104]{fk}, the $\Omega _i$ 
also form a basis of ${\cal H}^1(C)$. It follows that 
$\underline{\Omega}:z\mapsto\sum_{i=1}^3\, z_i\,\Omega _i$ 
is a linear isomorphism from $\C ^3$ onto ${\cal H}^1(C)$. 

The adjoint $\underline{\Omega}^*$ is a linear 
isomorphism from ${\cal H}^1(C)^*$ onto $\left(\C ^3\right) ^*\simeq \C^3$. 
For every standard basis vector $e_i$ of $\C ^3$ and closed curve 
$\gamma$ in $C$ we have that 
\[
\langle\underline{\Omega}^*\left(\oint_{\gamma}\right),\, e_i\rangle 
=\langle \oint_{\gamma},\,\underline{\Omega}\left( e_i\right)\rangle 
=\langle\oint_{\gamma},\, \Omega _i\rangle =\oint_{\gamma}\,\Omega _i,
\]
which shows that $\underline{\Omega}^*\circ\oint$ is equal to the mapping 
$\oint\Omega :\op{H}_1(C,\,\Z )\to\C ^3$. It follows that 
the image ${\cal P}\left(\Omega |_C\right)$ of  $\oint\Omega$ is equal to 
$\underline{\Omega}^*(\Lambda )$, which implies that ${\cal P}\left(\Omega |_C\right)$ 
is a lattice in $\C ^3$, because $\underline{\Omega}^*$ is a linear isomorphism and 
$\Lambda$ is a lattice in ${\cal H}^1(C)^*$. 

The second statement follows from the fact that the mapping 
\[
q\mapsto\left(\theta\mapsto\int_{q_0}^q\,\theta\right) : C\to {\cal H}^1(C)^*/\Lambda 
\]
is an embedding from $C$ into $\op{Jac}(C)$, cf. 
Farkas and Kra \cite[III.6.4]{fk}. 
\end{proof}

The pre\--image of $(\int\Omega )(C)$ in $\C ^3$ under the projection 
from $\C ^3$ to $\C ^3/{\cal P}\left(\Omega |_C\right)$ is a closed and smooth 
complex one\--dimensional analytic submanifold of $\C ^3$, 
which we denote by $\Gamma$. Note that $\Gamma$ is equal to the isotropic 
complex analytic curve which is associated to the minimal 
surface $M$ as in Section \ref{weierstrasssubsec}. 
The restriction to $\Gamma$ of the projection $\op{Re}:\C ^3\to\R ^3$ 
is a smooth immersion from $\Gamma$ onto the minimal surface $M$ in 
$\R ^3$ which is obtained by means of integration of $\omega$ 
over $C$. 

We also have that $\op{Re}{\cal P}\left(\Omega |_C\right) ={\cal P}$, 
the group of periods of $M$. Because 
${\cal P}\left(\Omega |_C\right)$ contains an 
$\R$\--basis of  $\C ^3$ and the real linear mapping $\op{Re}$ is 
surjective, it follows that ${\cal P}$ contains an $\R$\--basis of 
$\R ^3$. However, it can easily happen that ${\cal P}$ is 
not closed in $\R ^3$, in which case the minimal surface 
$M$ is dense in $\R ^3$. 

The period group ${\cal P}$ is closed, i.e. a lattice, in $\R ^3$, 
if and only if  the kernel $\left(\op{i}\R ^3\right)$ of 
the projection $\op{Re}:\C ^3\to\R ^3$ 
has an $\R$\--basis which consists of elements of ${\cal P}\left(\Omega |_C\right)$. 
In other words, if and only if ${\cal P}\left(\Omega |_C\right)$ 
contains three linearly independent 
purely imaginary elements. In this case the minimal surface 
$M$ is properly immersed in $\R ^3$, and is triply periodic. 

Extend the mapping $\oint\Omega :\op{H}_1(C,\,\Z )\to \C ^3$ 
to an $\R$\--linear mapping from 
$\oint\Omega :\op{H}_1(C,\,\R )\to \C ^3$, 
which we denote with the same letter. Then 
${\cal P}={\cal P}(C)$ is a lattice in $\R ^3$, 
if and only if the 3\--dimensional 
real linear subspace 
\[
L(C)=(\oint\Omega )^{-1}\left(\op{i}\R ^3\right)
\]
of the 6\--dimensional real vector space $\op{H}_1(C,\,\R )$ 
contains three linearly independent elements of the lattice 
$\op{H}^1(C,\,\Z )$. Let $G$ be the Grassmann manifold 
of all 3\--dimensional real linear subspaces of $\op{H}_1(C,\,\R )$,
and let $G_0$ denote the set of all $L\in G$ such that 
$L\cap \op{H}_1(C,\,\Z )$ contains a basis of $L$. 
Then $G_0$ is dense in $G$. Because of the discrete 
nature of the homology groups with values in $\Z$, 
there is a canonical identification, for any hyperelliptic curve 
$C'$ near $C$, of $\op{H}_1(C',\,\Z)$ with $\op{H}_1(C,\,\Z )$, 
and therefore also of $\op{H}_1(C',\,\R )$ with $\op{H}_1(C,\,\R )$, 
and of the Grassmann manifold of $\op{H}_1(C',\,\R )$ with $G$. 

Let ${\cal C}_0$ denote the set of $C'\in {\cal C}$ such that 
${\cal P}(C')$ is a lattice in $\R ^3$. If the mapping $C'\mapsto L(C')$ 
from the real 18\--dimensional space ${\cal C}$ to the real 9\--dimensional 
manifold $G$ has surjective derivative at $C$, 
then the conclusion is that ${\cal C}_0$ 
is dense in a neighborhood of $C$ in ${\cal C}$. 
Moreover, near $C$ the set ${\cal C}_0$  then is equal 
to the union of countably many smooth real analytic submanifolds of 
${\cal C}$, each of real dimension $18-9=9$. Of these 9 dimensions, 5 are 
essential, if we substract the 4 dimensions of the real rotations and the 
dilations. 

Pirola \cite{pirola} has proved that for the generic 
hyperelliptic curves $C$ as above the 
derivative of the mapping $C'\to L(C')$ surjective at $C$. 
It follows that the set of $C$, such that the period group 
${\cal P}$ is a lattice in $\R ^3$, is dense in ${\cal C}$. 
Moreover, in the open dense subset where $\op{D}\! L(C)$ is 
surjective, it is equal to the union of countably many smooth 
manifolds as described above. 

\medskip\noindent
There are many examples known of triply periodic minimal surfaces,  
among which those constructed by H.A. Schwarz and Alan Schoen, 
cf. DHKW \cite[pp. 212-217]{dhkw}. 
These minimal surfaces are not only properly immersed, 
but even embedded in $\R ^3$. 

The embedded minimal surfaces in the 
real five\--dimensional family of Meeks \cite[Thm. 7.1]{meeks} 
correspond to the genus three hyperelliptic curves $C$ as 
above, which satisfy the additonal condition that $C$ is 
invariant under the anti\--holomorphic involution
\[
A:(u,\, q)\mapsto\left( -1/\overline{u},\, 
\overline{q}/\overline{u}^{4}\right) 
\]
of $Q$ which is induced by the antipodal mapping $s\mapsto -s$ 
on the sphere $S$.  
The $[\gamma]\in\op{H}_1(C,\,\Z )$ such that 
$A ([\gamma ])=[\gamma ]$ form a three\--dimensional sublattice 
$L$ of $\op{H}_1(C,\,\Z )$. Furthermore, 
because $A^*\Omega =-\overline{\Omega}$, it follows that 
$(\oint\Omega )(L)\subset 
\op{i}\R ^3$, which implies that ${\cal P}$ is a lattice in $\R ^3$. 
The condition for $C$ means that 
\[
f_0(u)=c\, \prod_{j=1}^4\,\left( u-a_j\right)\, \left( u+1/\overline{a_j}\right) ,
\]
where $c\in\C$, $a_j\in\C$, $c\,\prod_{j=1}^4\, a_j\in\R$, and 
$a_j\neq a_k$, $a_j\neq -1/\overline{a_k}$ when $j\neq k$. 
This condition implies that the branch points consist of  
four arbitrary pairs of antipodal points on $S$, and there 
is a free nonzero real factor in $f_0$. 

For many of the known examples, 
Karcher and Wohlgemut in \cite[pp. 317-347]{karcher} 
found explicit Weierstrass data. With the substitutions 
$g=u$, $\frac{\op{d}\! h}{g}=\frac{-2\op{d}\! u}{q}$, 
$\op{d}\! h=\mu\,\frac{\op{d}\! g}{g}$, 
cf. Karcher \cite[p. 315 and 318]{karcher}, which implies that 
$\mu =-u^2/q$, we find that the examples in Karcher 
\cite[pp. 319-329]{karcher} belong to the family of 
Meeks \cite[Thm. 7.1]{meeks}. (In the case of (H) in 
Karcher \cite[p.329]{karcher}, the quantity 
$\mu ^2$ probably should be replaced by $\mu ^{-2}$.)  
The first example of Schwarz corresponds  
to the hyperelliptic curve $q^2/4=u^8-14\, u^2+1$, 
cf. Remark \ref{FWrem} and DHKW 
\cite[p. 175]{dhkw}. On the other hand, the examples in Karcher 
\cite[pp. 330-347]{karcher} do not look like 
hyperelliptic curves of genus three.

\subsection{The Costa Surface}
The famous {\em Costa surface} is a nonperiodic, embedded minimal surface 
of genus one and with three flat points at infinity. Furthermore the 
degree of its Gauss map is equal to three, cf. DHKW \cite[p. 198]{dhkw}. 

An explicit Weierstrass representation has been given 
in Hoffman and Meeks \cite[(1.1), (3.1), (3.5)]{hm}. If, in the case $k=1$, 
we substitute $g=u$, $w=c/u$, $\eta =-2\op{d}\! u/q$ and eliminate 
the variable $z$, then we arrive at the equation 
\begin{equation}
c^5\, q^3-c^2\, u^4\, q^2-27\, c^4\, u^4+4u^8=0
\label{Costa}
\end{equation}
for the curve $C$ in $(u,\, q)$\--coordinates. 
Therefore the degree of the Gauss map is equal to three. 
With the substituion $u=v^{-1}$, 
$q=v^{-4}\, r$, this equation is equivalent to 
\[
c^5\, r^3-c^2\, r^2-27\, c^4\, v^8+4v^4=0.
\]
According to Hoffman and Meeks \cite[prop. 3.1]{hm}, there 
is a unique positive value of $c$ such that the 
corresponding minimal surface has no periods. 

At $u=q=0$, which lies over $e_3\in S$, 
$C$ has intersection number with the fiber and the 
zero section equal to three and four respectively, and therefore 
$u=q=0$ corresponds to a flat point at infinity of $M$. At $v=r=0$, 
which lies over $-e_3\in S$, $C$ has two components, each of which 
has intersection number with the fiber and the zero section 
equal to one and two, respectively. Therefore $v=r=0$ corresponds to 
two flat points at infinity of $M$. The other intersection points 
with the zero section correspond to the equations 
$q=0$, $u^4=27\, c^4$. At each of these four points the intersection 
number with the fiber and the zero section is equal to two and one, 
respectively, which means that these points correspond to the four 
finite flat points of $M$. 

According to (\ref{genusspin}), the genus of $D$ is equal to 
$g=1+3-\frac{3}{2}-\frac{3}{2}=1$. In this way we have recovered the 
properties mentioned in the beginning of this section from the equation 
(\ref{Costa}) which defines $C$.

\medskip\noindent
{\bf Acknowledgement}~~~I thank Luc Vrancken, Joop Kolk, Erik van den Ban, 
Eduard Looijenga and Frans Oort for stimulating discussions.

\noindent{\bf Adress of the author}\\
Department of Mathematics, Utrecht University\\
Postbus 80 010, 3508 TA Utrecht, The Netherlands\\
e-mail: duis@math.uu.nl


\begin{thebibliography}{99}
\bibitem{bpv}W. Barth, C. Peters and A. Van de Ven: 
{\em Compact Complex Surfaces}. 
Springer\--Verlag, Berlin, etc., 1984.

\bibitem{c}C. Chevalley: {\em Theory of Lie Groups, Vol. I}. 
Princeton University Press, Princeton, New Jersey, 1946. 

\bibitem{conway}J.B. Conway: {\em Functions of One Complex 
Variable}. Springer\--Verlag, Berlin, Heidelberg, 1978. 

\bibitem{dhkw}U. Dierkes, S. Hildebrandt, A. K\"uster 
and O. Wohlrab: {\em Minimal Surfaces I}. Springer\--Verlag, 
Berlin, Heidelberg, 1992

\bibitem{dk}J.J. Duistermaat and J.A.C. Kolk: 
{\em Lie Groups}. Springer\--Verlag, Berlin, Heidelberg, 2000. 

\bibitem{fk}H.M. Farkas and I. Kra: {\em Riemann Surfaces}.  
Springer\--Verlag, New York, etc., 1992. 

\bibitem{gack}F. Gackstatter: {\em \"Uber die Dimension einer Minimalfl\"ache 
und zur Ungleichung von St. Cohn\--Vossen}. 
Archive for Rational Mechanics and Analysis  
{\bf 61} (1976) 141-152. 

\bibitem{gh}P. Griffiths and J. Harris: {\em Principles of 
Algebraic Geometry}. John Wiley \& Sons, New York, etc., 1978. 

\bibitem{h}F. Hirzebruch: {\em \"Uber eine Klasse von 
einfach\--zusammenh\"angenden komplexen Mannigfaltigkeiten}.
Mathematische Annalen {\bf 124} (1951) 77-86. 

\bibitem{hk}D. Hoffman and H. Karcher: 
{\em Complete Embedded Minimal 
Surfaces of Finite Total Curvature}. pp. 5-93 in 
R. Osserman (Ed.): {\em Geometry V}. Encyclopaedia of Mathematical 
Sciences, Volume 90. Springer\--Verlag, Berlin, Heidelberg, 1997. 

\bibitem{hm}D. Hoffman and W.H. Meeks III: 
{\em Embedded minimal surfaces of finite topology}. 
Annals of Mathematics {\bf 131} (1990) 1-34.

\bibitem{jm}L. Jorge and W.H. Meeks III: {\em The topology of 
complete minimal surfaces of finite total Gaussian curvature}. 
Topology {\bf 22} (1983) 203-221. 

\bibitem{karcher}H. Karcher: {\em The triply periodic minimal 
surfaces of Alan Schoen and their constant mean curvature companions}. 
Manuscripta mathematica {\bf 64} (1989) 291-357. 

%\bibitem{lawson}H.B. Lawson Jr.: {\em Lectures on Minimal Submanifolds}. 
%Publish or Perish, Berkeley, 1971.
%

\bibitem{meeks}W.H. Meeks III: {\em The theory of triply 
periodic minimal surfaces}. Indiana University Mathematics Journal 
{\bf 39} (1990) 877-936. 

\bibitem{l}S. {\L}ojasiewicz: {\em Introduction to Complex 
Analytic Geometry}. Birkh\"auser Verlag, Basel, Boston, Berlin, 1991. 

\bibitem{o59}R. Osserman: {\em An analogue of the Heinz\--Hopf inequality}. 
J. Math. Mech. {\bf 8} (1959) 383-385. 

\bibitem{o63}R. Osserman: {\em On complete minimal surfaces}. 
Archive for Rational Mechanics and Analysis {\bf 13} (1965) 392-404. 

\bibitem{o}R. Osserman: {\em A Survey of Minimal Surfaces}. 
Van Nostrand Reinhold Company, New York, etc., 1969. 

\bibitem{pirola}G.P. Pirola: {\em The infinitesimal variation of the 
spin abelian differentials and periodic minimal surfaces.} 
Comm. Anal. Geom. {\bf 6} (1998) 393-426.  

%\bibitem{s}R. Schoen: Uniqueness, symmetry, and embeddedness of 
%minimal surfaces. {\em Journal of Differential Geometry} 
%{\bf 18} (1983) 791-809. 

\bibitem{sha}I.R. Shafarevich: {\em Basic Algebraic Geometry}. 
Springer\--Verlag, Berlin, Heidelberg, New York, 1977. 
\end{thebibliography}
\end{document}